\documentclass[12pt]{amsart}
\usepackage{graphicx}
\usepackage{amsmath}
\usepackage{amsfonts}
\usepackage{amssymb}
\usepackage{amscd}

\newcommand{\la}{\lambda}
\newcommand{\be}{\beta}
\newcommand{\ga}{\gamma}
\newcommand{\de}{\delta}
\newcommand{\De}{\Delta}
\newcommand{\al}{\alpha}
\newcommand{\e}{\varepsilon}

\newcommand{\si}{\sigma}
\newcommand{\Si}{\Sigma}
\newcommand{\BR}{\mathbb{R}}
\newcommand{\BZ}{\mathbb{Z}}
\newcommand{\BQ}{\mathbb{Q}}
\newcommand{\BN}{\mathbb{N}}
\newcommand{\A}{\mathcal{A}}
\newcommand{\U}{\mathcal{U}}
\newcommand{\m}{\mathfrak{m}}
\newcommand{\ov}{\overline}

\newcommand{\X}{\mathfrak{X}}
\newcommand{\R}{\mathfrak{R}}
\renewcommand{\P}{\mathcal{P}}

\newcommand{\Ral}{\mathcal{R}_\al}
\newcommand{\Xal}{\mathcal{X}_\al}
\newcommand{\T}{\mathcal{T}_\al}
\newcommand{\V}{\mathcal{V}_\al}
\newcommand{\Leb}{\mathcal{L}}

\newcommand{\wt}{\widetilde}
\newcommand{\0}{\mathbf{0}}
\renewcommand{\t}{\mathbf{t}}
\newcommand{\Tr}{\mathrm{Tr}}

\newcommand{\BT}{\mathbb{T}}

\renewcommand{\A}{\mathcal{A}}
\renewcommand{\phi}{\varphi}
\renewcommand{\c}{\mathfrak{c}}

\newcommand{\pic}[4]{
\begin{figure}[h]
\begin{center}
\includegraphics[origin=br,width=#1]{#2}
\end{center}
\caption{#3} \label{#4}
\end{figure}}

\newtheorem{lemma}{Lemma}[section]

\newtheorem{prop}[lemma]{Proposition}
\newtheorem{thm}[lemma]{Theorem}
\newtheorem{cor}[lemma]{Corollary}
\theoremstyle{definition}
\newtheorem{Def}[lemma]{Definition}
\newtheorem{exam}[lemma]{Example}
\newtheorem{conjecture}[lemma]{Conjecture}

\newtheorem{question}[lemma]{Question}
\theoremstyle{remark}
\newtheorem{rmk}[lemma]{Remark}

\begin{document}

\title[Arithmetic Dynamics]{Arithmetic Dynamics}
\author{Nikita Sidorov}
\address{Department of Mathematics, UMIST, P.O. Box 88,
Manchester M60 1QD, United Kingdom. E-mail:
Nikita.A.Sidorov@umist.ac.uk}
\date{\today}
\thanks{Supported by the EPSRC grant no GR/R61451/01.}
\dedicatory{To Paul Glendinning on the occasion of the birth of
his twins}

\begin{abstract} This survey paper is aimed to describe a
relatively new branch of symbolic dynamics which we call
Arithmetic Dynamics. It deals with explicit arithmetic
expansions of reals and vectors that have a ``dynamical" sense.
This means precisely that they (semi-) conjugate a given
continuous (or measure-preserving) dynamical system and a
symbolic one. The classes of dynamical systems and their codings
considered in the paper involve:
\begin{itemize}
\item Beta-expansions, i.e., the radix expansions in non-integer
bases;
\item ``Rotational" expansions which arise in the problem
of encoding of irrational rotations of the circle;
\item Toral expansions which naturally appear in arithmetic symbolic
codings of algebraic toral automorphisms (mostly hyperbolic).
\end{itemize}
We study ergodic-theoretic and probabilistic properties of these
expansions and their applications. Besides, in some cases we
create ``redundant" representations (those whose space of
``digits" is {\em a priori} larger than necessary) and study
their combinatorics.
\end{abstract}

\maketitle

\footnotesize
\tableofcontents
\normalsize

\section{Introduction}

The present survey paper is devoted to the recent progress in
the new branch of dynamical systems theory which we call {\em
Arithmetic Dynamics} (AD). It is worth pointing out that there
is no conventional agreement regarding what the expression AD
actually stands for -- at the moment when I am writing these
words (the year~2002), different people seem to see it
differently (see, {\em e.g.}, \cite{Miles, EEW}). This is not
particularly surprising: since the term is not fixed, anything
that has something to do with arithmetics and dynamics may
qualify.

Nonetheless, in the present paper the scope will be more narrow
than that; namely, we define AD as a discipline which deals with
symbolic codings of continuous (or measure-preserving) dynamical
systems (invertible or not) expressed in terms of {\bf explicit}
arithmetic expansions of real numbers of vectors. If one accepts
this (rather vague) definition, the classic ergodic theory of
continued fractions, for instance, quite fits into the scope of
AD. The term in question in this setting was suggested first by
A.~Vershik in the mid-1990's (implying that it would be suitable
for a full-developed theory in the future).

The expansions in questions will be usually called {\em
arithmetic codings}. Here are some characteristic features of
arithmetic codings:
\begin{itemize}
\item They extensively use number-theoretic methods and
techniques;
\item The arithmetic structure possessed by dynamical systems we
encode, is fully preserved;
\item As long as there are no obstacles of number-theoretic
nature, they can be generalized.
\end{itemize}

Apart from the ``normal" (i.e., one-to-one a.e.) expansions, there
exists another important topic that may be regarded as a part of
AD, namely, the theory of {\em redundant} or {\em excessive}
representations of real numbers and vectors. The model is as
follows: assume we have a fixed ``normal" expansion and the
natural set of its ``digits" is not a Cartesian product but has,
for instance, pairwise (or even more sophisticated) restrictions.
Our goal is to study the pattern, in which we lift all the
restrictions between {\bf distinct} digits and leave only the
minimal Cartesian hull for the original set of digits. For
example, in the case of $\be$-expansions with a non-integer $\be$
(see below) this eventually leads to the well-known {\em Bernoulli
convolutions}. Special attention will be paid to the combinatorics
of {\bf all} representations of a given $x$ as well as to the set
of those $x$ which despite lifting the restrictions, will have a
{\em unique representation} in the class in question.

One important point has to be made: AD is still in the cradle,
so to speak, i.e., definitely not yet a ``full-scale" subarea of
Dynamical Systems. Here are its lacks:
\begin{itemize}
\item At present there is hardly any systematic approach that
would cover a more or less substantial variety of
measure-preserving (or continuous) transformations. On the
contrary, as we will see, for each class of maps under
consideration the model turns out to be state-of-the-art.
\item Another serious issue regarding AD is that the
constructions in question are not yet quite robust and strongly
depend on the arithmetic structure of a dynamical system in
question.
\end{itemize}

Despite all this, most constructions look rather nice and are
closely related to number-theoretic problems as well (especially
for arithmetic codings of toral automorphisms). Thus, we believe
that even in this intermediate state AD is worth a detailed
description, with the hope that some day it will become an area
with a more systematic approach. Such a description is the aim
of the present paper.

Our intention is mostly to summarize the progress in AD in the
recent 10 years. The reason for this particular figure is
A.~Vershik's seminal paper \cite{Ver-ML} which appeared in the
winter of 1991--92 and which has stimulated quite a number of
new works and a great deal of ideas in AD. This paper presents
an arithmetic coding of the Fibonacci automorphism of the
2-torus (see Section~\ref{AC} for the definitions) based on the
two-sided generalization of the corresponding adic
transformation suggested by the same author in the late 1970's
\cite{Ver81, Ver82} (see Section~\ref{RE} and Appendix).

The structure of the paper is as follows:
\textbf{Section~\ref{BE}} is devoted to the {\em
$\be$-expansions}, i.e., radix representations in (generally
speaking, non-integer) bases $\be$ with $\be>1$. In particular,
we will briefly outline some well-known results about the {\em
greedy} and {\em lazy} beta-expansions, then describe in detail
the recent progress in the theory of {\em unique}
beta-representations in a fixed alphabet and finally, will deal
with the space of {\em all } beta-representations of a given
real number. Besides, we are going to devote a part of this
section to the case of the one-parameter family of {\em
intermediate} beta-expansions, i.e., those which lie in between
greedy and lazy ones (in the sense of the lexicographic
ordering).

\textbf{Section~\ref{RE}} deals with, generally speaking, quite
a different subject, namely, arithmetic codings of a dynamical
system with a purely discrete cyclic spectrum. More precisely,
we present two models for adic realization of an irrational
rotation of the circle $S^1=\BR/\BZ$. Both schemes deal with
expansions of the elements of $S^1$ in bases involving the
sequence of best approximation of the angle of rotation, which
we call {\em rotational expansions}. In this case the compacta
we obtain, are, generally speaking, non-stationary, and the map
is not a traditional shift but the {\em adic transformation}
(see Appendix for the definition).\footnote{There is nonetheless
some intersection with Section~\ref{BE} -- the two models
coincide if the angle of rotation is $\frac12(-a+\sqrt{a^2+4})$
for some $a\in\BN$.} This map is in a way transversal to the
shift (if the latter is well defined). Putting it simply, the
adic transformation is a generalization of the ``adding machine"
in the ring of $p$-adic integers to more general Markov
compacta.

We study the distribution of the ``digits" in both cases and
give sufficient conditions for the laws of large numbers and the
central limit theorem to hold for them. In the end of
Section~\ref{RE} we -- similarly to Section~\ref{BE} -- consider
the set of unique rotational expansions, prove a number of
claims that fully describe it and discuss the effect of
non-stationarity, which leads to an essential difference with
the beta-expansions.

Finally, in \textbf{Section~\ref{AC}} we consider arithmetic
codings of hyperbolic automorphisms of a torus; as was mentioned
above, it was a discovery in this area that initiated most of
the research in AD in the last 10 years. After the model case of
the {\em Fibonacci automorphism} of $\BT^2$, i.e., the one given
by the matrix $\begin{pmatrix} 1&1\\1&0\end{pmatrix}$, was
studied in detail in \cite{Ver-ML}, A.~Vershik \cite{Ver} asked
the question whether it is possible to generalize this
construction to all hyperbolic (or even ergodic) automorphisms.
More precisely, is it possible to find a symbolic coding of a
given automorphism $T$ such that certain structures (the stable
and unstable foliation, homoclinic points) have a clear
expression in the corresponding symbolic compactum. As is well
known, the classical models that deal with Markov partitions
\cite{Sinai, Bowen, GurSinai} do not have this property (apart
from possibly the case of dimension two).

Thus, in Section~\ref{AC} we describe the evolution of this area
in the past 10 years. This includes the Fibonacci automorphism
(as well as ergodic automorphisms of the 2-torus considered in
\cite{SV2, VS2}), the Pisot automorphisms, i.e. those whose
unstable (or stable) foliation is one-dimensional, and finally,
we will try to summarize all attempts to cope with the general
hyperbolic case and the difficulties that occur in doing so.

\textbf{Appendix} serves mostly auxiliary purposes: it briefly
describes the theory of adic transformations.

The experienced reader will notice, of course, that some topics
that might be in this survey paper, are missing. The reason for
this is that either a nice exposition of the corresponding
theory and results can be found elsewhere or the similarity with
the actual framework of AD, however broad it is, is imaginary.
Here is the list of some areas and subjects in question:
\begin{itemize}
\item {\em Odometers}. By this word is usually implied the adic
transformation (see Appendix) on a general compact set (not
necessarily a Markov compactum) with some extra ``fullness"
conditions -- see, {\em e.g.}, \cite{GLT}. I am still not
convinced though that the odometers are really natural: as far
as I am concerned, there are no interesting examples of symbolic
codings of non-symbolic dynamical systems by means of odometers.
\item {\em Algebraic codings of higher-rank actions on a torus}.
The reader may refer to the paper by Einsiedler and Schmidt
\cite{ES}, in which a general construction (similar to the one
from \cite{Sch} discussed in Section~\ref{AC}) was suggested.
There are also some partial results in this direction that can be
found in the author's paper \cite{S2}.
\item {\em Measures of arithmetic nature}. By those we mean
mainly Bernoulli-type convolutions parameterized by a Pisot
number. The information about them can be found, for example, in
\cite{OST} (see also references therein).
\end{itemize}

Since this is a survey paper, the proofs in the text are usually
omitted; exception is made for a few new results, namely,
Theorem~\ref{gencase} and all claims in Section~\ref{URE} as
well as some more minor claims.

The author is indebted to Anatoly Vershik for helpful
suggestions, remarks and historical data.

\section{Beta-expansions}\label{BE}

Let $\be>1$ and $x\ge0$; we call any representation of the form

\begin{equation}
x=\pi_\be(\e)=\sum_{n=1}^\infty \e_n\be^{-n} \label{beta-exp}
\end{equation}
a {\em $\be$-expansion} (or {\em beta-expansion} if we do not
want to specify the base). We need to make one historical
remark. The matter is that in the literature the beta-expansion
is often the one that uses the greedy algorithm for obtaining
the ``digits" $\e_n$ (see Section~\ref{greedy}); we acknowledge
this tradition (which apparently dates back to B.~Parry's
pioneering work \cite{Pa}), but think that within this text it
would only cause confusion, because we intend to describe
different types of algorithms for one and the same $\be$.

\subsection{Greedy expansions and the beta-shift}\label{greedy}
This subsection is one of the few exceptions we have made from
the rule that only the recent progress will be discussed. The
reason for doing so is simple: without clear exposition of this
model, it is impossible to explain the importance of ideas
developed in recent papers. In this subsection we will confine
ourselves mostly to the case of $x<1$ (although beta-expansions
can be extended to the positive half-axis as well -- see
Lemma~\ref{half}).

So, let $\tau_\be$ be the {\em $\be$-transformation}, i.e., the
map from $[0,1)$ onto itself acting by the formula
$$
\tau_\be(x)=\be x\bmod1
$$
(see Figure~\ref{Beta-shift}).

\pic{8cm}{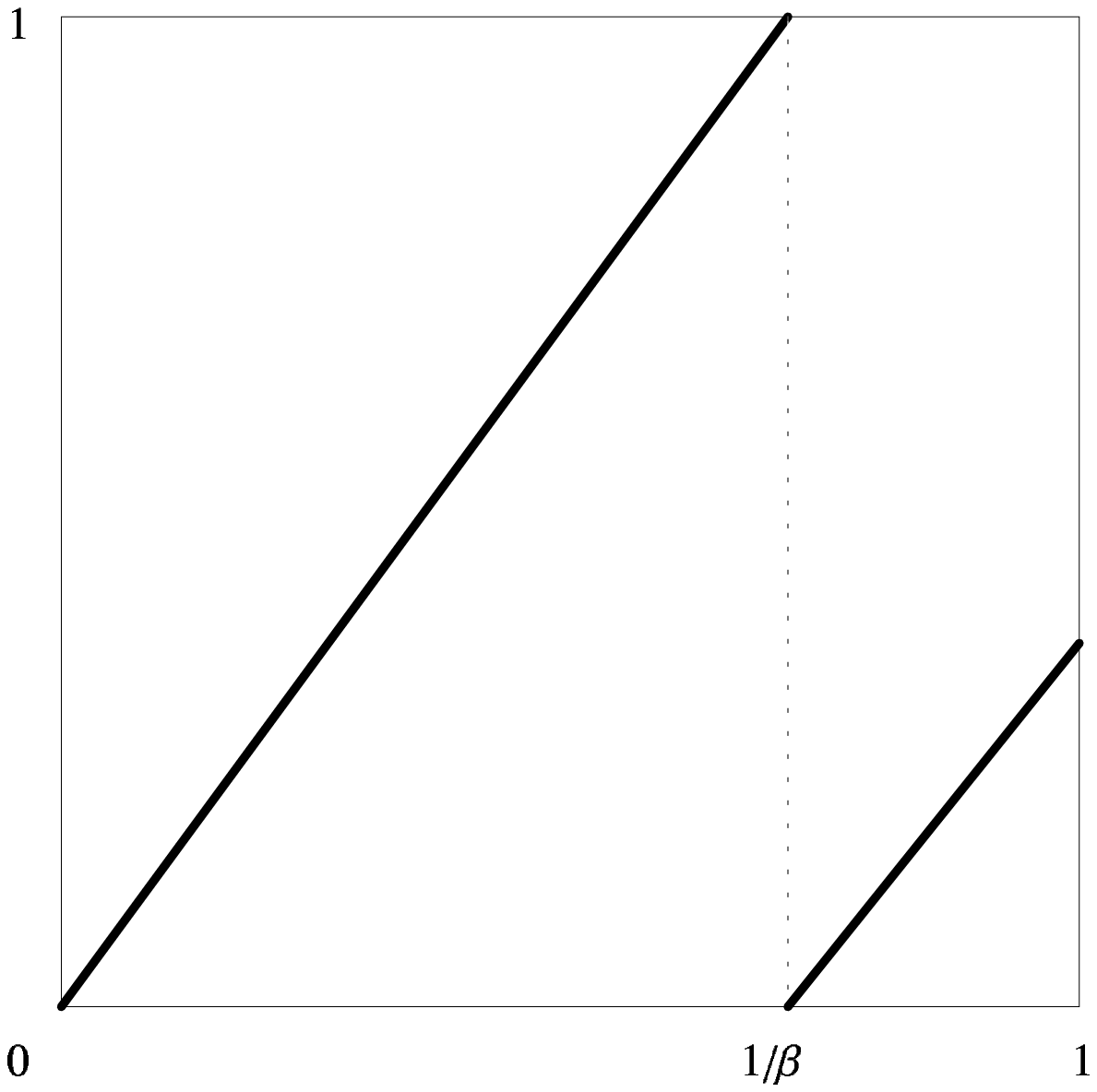}{The $\be$-shift}{Beta-shift}

This map is very important in ergodic theory and as well as
number theory, and the main tool for its study is its symbolic
encoding which we are going to describe below. Note that if
$\be=d$ is an integer, then $\tau_\be$ is isomorphic to the full
shift on $d$ symbols. The idea suggested in \cite{Re} was to
generalize this construction to the non-integer $\be$'s.

As is well known, to ``encode" it, one needs to apply the greedy
algorithm in order to obtain the digits in (\ref{beta-exp}),
namely, $\e_n=[\be\tau_\be^{n-1}x],\ n\ge1$. Then the one-sided
shift $\si_\be$ in the space $X_\be$ of all possible sequences
$\e$ that can be obtained this way, is clearly isomorphic to
$\tau_\be$, with the conjugating map given by (\ref{beta-exp}).
We will call $\si_\be$ the {\em $\be$-shift}. It is obviously a
proper subshift of the full shift on
$\prod_1^\infty\{0,1,\dots,[\be]\}$. The question is, what kind
of subshift is this -- or equivalently -- what is actually
$X_\be$?

This question was answered by B.~Parry in the seminal paper
\cite{Pa}. The theorem he proved is the following. Let the
sequence $(a_n)_1^\infty$ be defined as follows: let
$1=\sum_{1}^{\infty}a_k' \be^{-k}$ be the greedy expansion of 1,
i.e, $a_n'=[\be\tau_\be^{n-1}1],\ n\ge1$; if the tail of the
sequence $(a_n')$ differs from $0^\infty$, then we put
$a_n\equiv a_n'$. Otherwise let $k=\max\,\{j:a_j'>0\}$, and
$(a_1,a_2,\dots):= (a_1',\dots,a_{k-1}',a_k'-1)^\infty$.

\begin{thm} \cite{Pa}\label{Pabe}
\begin{enumerate}
\item For any sequence $(a_n)_1^\infty$ described above, each
power of its shift is lexicographically less than or equal to
the sequence itself, and the equality occurs if and only if it
is purely periodic. Conversely, each sequence having this
property is $(a_n(\be))_1^\infty$ for some $\be>1$.\label{P1}
\item For each greedy expansion $(\e_n)_1^\infty$ in base $\be,\
(\e_n,\e_{n+1},\dots)$ is lexicographically less (notation:
$\prec$) than $(a_1,a_2,\dots)$ for every $n\ge1$. Conversely,
every sequence with this property is actually the greedy
expansion in base $\be$ for some $x\in [0,1)$.
\end{enumerate}
\end{thm}

Let $\mathbf a=(a_n)_1^\infty$ and $\si$ stands for the general
one-sided shift on sequences. Thus, we have
\begin{equation}
\label{Xbe}
X_\be=\left\{\e\in\prod_1^\infty\{0,1,\dots,[\be]\}\mid
\si^n\e\prec\mathbf{a},\ n\ge0\right\}
\end{equation}
and the following diagram commutes:
\[
\begin{CD}
X_\be @>{\si_\be}>> X_\be \\
@V{\pi_\be}VV @VV{\pi_\be}V \\
[0,1) @>{\tau_\be}>> [0,1)
\end{CD}
\]
Note that usually $X_\be$ is called the (one-sided) {\em
$\beta$-compactum}.

Types of subshifts people know well how to deal with are mostly
SFT (subshifts of finite type) and their factors, called {\em
sofic subshifts} -- see, {\em e.g.}, \cite{LindMarcus}. It is
thus natural to ask whether $X_\be$ is such and if so, for which
$\be$? The following theorem gives a rather disappointing answer
to this question. Roughly speaking, only certain algebraic $\be$
yield sofic beta-compacta. Recall that an algebraic integer is
called a {\em Pisot number} if it is a real number greater than
1 and all its conjugates are less than 1 in modulus. A {\em
Perron number} is an algebraic integer $\be$ greater than 1
whose conjugates are less than $\be$ in modulus.

\begin{thm} \cite{Pa, Ber1, Lind} (see also the survey \cite{Bl})
\label{sofic}
\begin{enumerate}
\item The $\be$-compactum $X_\be$ is an SFT if and only if the
sequence $a'_n$ (see above) is finite (i.e., its tail is
$0^\infty$).
\item $X_\be$ is sofic if and only if $(a_n)$ is ultimately
periodic.
\item If $\be$ is a Pisot number, then $X_\be$ is sofic.
\item If $X_\be$ is sofic, then $\be$ is a Perron number.
\end{enumerate}
\end{thm}

It has to be said that a more or less explicit description of
$X_\be$ in case of transcendental $\be$ (as well as for $\be=3/2$,
say) seems to be hopeless -- see \cite{Bl}. However, the
ergodic-theoretic properties of the beta-shift are well studied
and clearly understood by now for all $\be>1$. The following
statement summarizes them.

\begin{thm} \cite{Re, Pa, Sm, Hof}
\begin{enumerate}
\item The $\be$-transformation is topologically mixing, and its
topological entropy is equal to $\log\be$.
\item The $\be$-shift is intrinsically ergodic\footnote{This means
by definition that the map has a unique measure of maximal
entropy.} for any $\be>1$. \label{ier}
\item The unique measure of maximal entropy for $\tau_\be$ is
equivalent to the Lebesgue measure on $[0,1]$ and the
corresponding density is bounded from both sides.
\item The natural extension of $\tau_\be$ is Bernoulli. Moreover,
the $\be$-shift $\si_\be$ is weakly Bernoulli with respect to
the natural (coordinate-wise) partition.
\end{enumerate}
\end{thm}

\begin{rmk} The proof of item~(\ref{ier}) given in \cite{Hof} is
rather complicated. In fact, I do not really understand why the
proof of the intrinsic ergodicity of the transitive subshifts of
finite type given by B.~Parry (see, {\em e.g.},
\cite[pp.~194--196]{Wal}) cannot be applied to the $\be$-shifts
as well. The only property one needs apart from ergodicity, is
the fact that the measure of any cylinder of length~$n$ divided
by $\be^n$ is uniformly bounded. This is well known since
\cite{Re}.
\end{rmk}
\begin{rmk} In \cite{Ber2} a simplified proof of the above theorem
was given. It is based on some auxiliary results on {\em coded
systems} (which the $\be$-shift is -- see, {\em e.g.},
\cite{Bl}). Unfortunately, this manuscript, rather helpful from
the methodological point of view, is unpublished and not very
easy to get hold of.\footnote{In our days if a manuscript is
unpublished, this is not necessarily that bad: a \TeX\ file is
even simpler to deal with than a hard copy. Unfortunately, this
particular manuscript dates back the pre-\TeX\ epoque\dots}
\end{rmk}

The greedy expansion can be alternatively characterized as
follows: assume that $n\ge2$ and that the first $n-1$ digits of
the expansion~(\ref{beta-exp}) are already chosen. Then if there
is a choice for $\e_n$, we choose the largest possible number
between 0 and $[\be]$. Similarly, if we choose the smallest
possible $\e_n$ every time when we have a choice, this expansion
is called the {\em lazy $\be$-expansion}. Let us formally explain
what we mean by the existence of a choice. Let
$r_n(x,\be):=x-\sum_1^{n-1}\e_k\be^{-k}$; if
$r_n(x,\be)<\be^{-n}$, then $\e_n$ has be to equal to 0. If, on
the contrary, $r_n(x,\be)\ge[\be]\be^{-n}$, then inevitably
$\e_n=[\be]$. As is easy to see, in any other case there will be a
choice for $\e_n$.

The following assertion is straightforward:

\begin{lemma} For a given $x\in[0,1)$, each of its
$\be$-expansions of the form~(\ref{beta-exp}) lies between its
lazy and greedy $\be$-expansions in the sense of lexicographic
ordering of sequences.\label{int}
\end{lemma}

We will return to the case of ``intermediate" expansions (i.e.,
those that lie strictly between the lazy and greedy ones) in
Section~\ref{intermediate}.

\begin{rmk}\label{nonn} As we have mentioned above, it is
possible to expand any positive number (not necessarily from
$(0,1)$) in base $\be$ by means of the greedy algorithm. Namely,
let $\mathbf a$ be as above, and
\begin{equation}\label{Xbetilde}
\wt
X_\be:=\left\{\e\in\prod_{-\infty}^\infty\{0,1,\dots,[\be]\}\mid
(\e_n,\e_{n+1},\dots)\prec\mathbf{a},\ n\in\BZ\right\},
\end{equation}
i.e., the natural extension of the $\be$-compactum. We will call
it the {\em two-sided $\be$-compactum}; it will be used
extensively in Section~\ref{AC}.
\begin{lemma}\label{half}
Any $x\ge0$ has the greedy $\be$-expansion of the form
\begin{equation}\label{32-sided}
x=\sum_{n=-\infty}^{+\infty} \e_n\be^{-n},
\end{equation}
where $(\e_n)$ is a sequence from $\wt X_\be$ finite to the
left, i.e., $\e_n\equiv0$ for $n\le N_0$ for some $N_0\in\BZ$.
\end{lemma}
The proof of this claim is similar to the one on the $p$-adic
representations of positive reals, and we omit it (see also
Section~\ref{AC}).
\end{rmk}

\subsection{Unique expansions and maps with gaps}\label{un} This
subsection is aimed to describe a branch of a relatively new
direction in the theory of arithmetic expansions, which deals with
lifting all restrictions on ``digits" leaving only the ``Cartesian
hull" (see Introduction).
\subsubsection{$1<\be<2$} Assume first that
$\be\in (1,2)$ and let $\Si=\prod_1^\infty\{0,1\}$. In the
previous subsection we have described the specific (greedy)
algorithm for choosing ``digits" $\e_n$. As we have seen, the
cost for this (very natural) approach is that the set of digits
is quite complicated and unless $\be$ is an algebraic number,
there is hardly any hope to describe it more or less explicitly
(see Theorem~\ref{sofic}).

In the 1990's a group of Hungarian mathematicians led by Paul
Erd\"os began to investigate 0-1 sequences that provide {\bf
unique} representations of reals \cite{ErdJoo, ErdJooKom,
ErdJooKom2}. More precisely, let
\[
\A'_\be=\left\{x\in\left(0,\frac1{\be-1}\right)\mid \exists !\
(\e_n)_1^\infty\in\Si : x=\sum_{n=1}^\infty\e_n\be^{-n}\right\}
\]
(it is obvious that the only representation for $x=0$ is
$0^\infty$ and the only representation for $x=1/(\be-1)$ is
$1^\infty$, so we will exclude both ends of the interval).

The first result about this set is given in \cite{ErdJoo}:
\begin{prop}
\label{Leb=0} The set $\A'_\be$ has Lebesgue measure zero for any
$\be\in(1,2)$. Moreover, if $\be<G$, where $G=\frac{1+\sqrt5}2$,
then in fact every $x$ has $2^{\aleph_0}$ representations in the
form~(\ref{beta-exp}).
\end{prop}

The question is, what can one say about the cardinality and --
in case it is the continuum --  about the Hausdorff dimension of
this set. The answer to this question is given by
P.~Glen\-din\-ning and the author in \cite{GS1}.

To present this result, we need some preliminaries. Let $\be_*$
denote the \textit{Komornik-Loreti constant} introduced by
V.~Komornik and P.~Loreti in \cite{KL}, which is defined as the
unique solution of the equation
$$
\sum_{1}^{\infty}\mathfrak{m}_{n}x^{-n+1}=1,
$$
where $\mathfrak{m}=(\mathfrak{m}_n)_1^\infty$ is the Thue-Morse
sequence
$$
\mathfrak{m}=0110\,\,1001\,\,1001\,\,0110\,\,1001\,\,0110\dots,
$$
i.e., the fixed point of the substitution $0\to01,\ 1\to10$. The
Komornik-Loreti constant is known to be transcendental
\cite{AllCos}, and its numerical value is approximately as
follows:
\[
\be_*=1.787231650\dots
\]
The reason why this constant was introduced in \cite{KL} is that
it proves to be the smallest number $\be$ such that $x=1$ has a
unique representation in the form~(\ref{beta-exp}). Now we are
ready to formulate the result we mentioned above.

\begin{thm}\cite{GS1} The set $\mathcal{A}'_\be$ is:
\begin{itemize}
\item empty if $\be\in (1,G]$;
\item countable for $\be\in(G,\be_*)$;
\item an uncountable Cantor set of zero Hausdorff dimension if
$\be=\be_*$; and
\item a Cantor set of positive Hausdorff dimension for
$\be\in (\be_*,2)$.
\end{itemize}
\end{thm}

The proof of this result given in \cite{GS1} is based on the
observation that for $\e$ to be a unique expansion for some
$x\in\left(\frac{2-\be}{\be-1},1\right)$, it has to be its both
greedy and lazy expansion -- this is a direct consequence of
Lemma~\ref{int}. Let
$\A_\be=\A'_\be\cap\left(\frac{2-\be}{\be-1},1\right)$ and
$\U_\be=\pi_\be^{-1}(\A_\be)$, where $\pi_\be$ is given by
(\ref{beta-exp}). It suffices to use (\ref{Xbe}) and a similar
condition for the lazy expansion, which leads to the following
lemma on the structure of the set $\U_\be$:

\begin{lemma}\cite{GS1} \label{Ubeta} The set $\U_\be$ can be
described as follows:
\begin{equation}
\U_\be=\{\e\in\Si:\ov{\mathbf{a}}\prec\si^n\e\prec\mathbf{a},\
n\ge0\}. \label{Uq}
\end{equation}
\end{lemma}

\begin{rmk} It is worth noting that although (\ref{Xbe}) looks
similar to (\ref{Uq}), the compacta $X_\be$ and $\U_\be$ are
completely different. In particular, the entropy of the
$\be$-shift is $\log\be$ and in the case of the shift on
$\U_\be$ it is constant a.e. -- see Theorem~\ref{erg-prop}
below. Note also that restrictions like (\ref{Uq}) are quite
common in one-dimensional dynamics, and this is not a
coincidence -- see \cite{GS2}.
\end{rmk}

From Lemma~\ref{Ubeta} one can obtain the result on the
cardinality of $\U_\be$ and therefore of $\A_\be'$ as well --
see \cite{GS1}.

The problem with this proof is that it may be characterized as
``a rabbit out of a hat" type of proof. Indeed, it does not say
anything about the origin of the Komornik-Loreti constant as the
main threshold between countable and uncountable set of uniquely
representable points. This collision was overcome in the
subsequent papers by the same authors \cite{GS2, GS3}.

The key paper \cite{GS2} is devoted to a ``dynamical" version of
the proof which will also involve various ergodic-theoretic and
geometric applications for the shift on the space $\U_\be$. Let
us explain, where dynamics enters the game.

In a number of recent works {\em maps with gaps} or a {\em maps
with holes} have been considered -- see, {\em e.g.}, \cite{BKT}
and references therein. The model in question is as follows: let
$X\subset\BR^d$ and $f:X\to X$ be a map (invertible or not) with
positive topological entropy. Let $D$ be an open subset of $X$.
The idea is to study the ``dynamics of $f$ on $X\setminus D$".
More precisely, let
\[
U=X\setminus \bigcup_{n\in\BZ}f^n(D)\quad\textrm{if $f$ is
invertible}
\]
or
\[
U=X\setminus \bigcup_{n\le0}f^n(D)\quad\textrm{if $f$ is
non-invertible}.
\]
One may ask two questions about this model.

\medskip\noindent {\bf Question 1.} Is $U$ empty? countable?
uncountable? a Cantor set of po\-sitive Hausdorff dimension?

\medskip\noindent {\bf Question 2.} If $U$ has positive
Hausdorff dimension, describe the dynamics of $f_U=f|_U$
(sometimes called the {\it exclusion map}). For instance, is
$f_U$ transitive? topologically mixing? intrinsically ergodic?
etc.

Some (mostly, ``generic" results) in this direction for Axiom~A
maps $f$ on smooth manifolds can be found in \cite{CM1, CM2,
CMT, BKT}. In particular, in \cite{BKT} it is shown that if $f$
is a hyperbolic algebraic automorphism of the torus $\BT^m$ (see
Section~\ref{AC} for the relevant definitions) and
$D=D(a_1,\dots, a_m)$ is the parallelepiped built along the
leaves of the stable and unstable foliation of $f$ passing
through $\0$ with the sides of length~$a_1,\dots,a_m$, then for
a Lebesgue-generic $m$-tuple $(a_1,\dots,a_m)$ the exclusion map
$f|_U$ is a subshift of finite type.

Return to our situation. Let $T_\be$ be the following ``map with
a gap" acting from $[0,1/(\be-1)]$ onto itself:
\begin{equation}
T_\be(x)=\begin{cases} \be x,& x\in\bigl[0,\frac1{\be}\bigr]\\
                       \textrm{not defined},& x\in
                  \bigl(\frac1{\be},\frac1{\be(\be-1)}\bigr)\\
                       \be x-1,& x\in
                  \bigl[\frac1{\be(\be-1)},\frac1{\be-1}\bigr]
\end{cases}
\label{Tbetaa}
\end{equation}
(see Figure~\ref{Mapgap}).

\pic{9truecm}{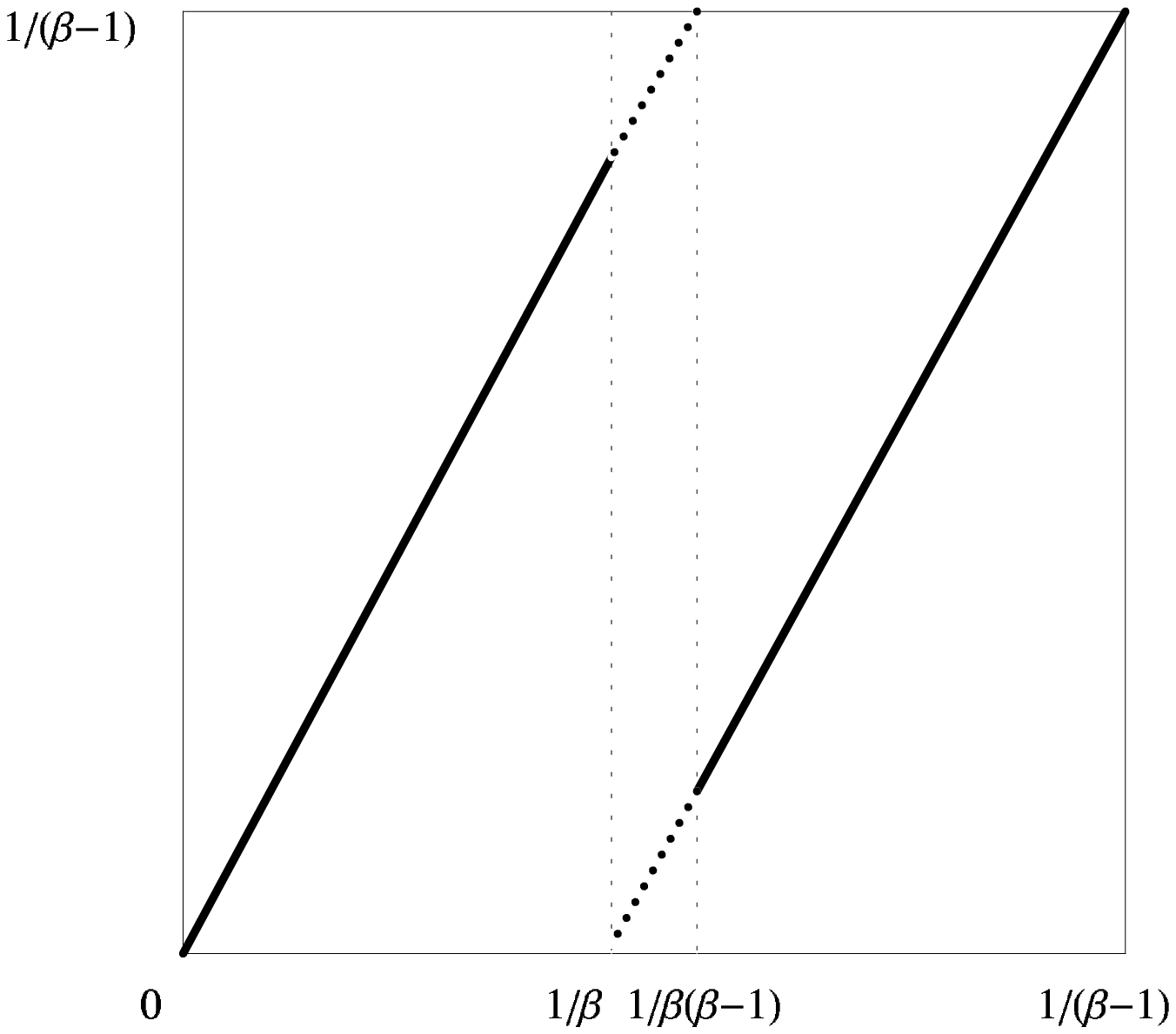}{The map with a gap $T_\be$}{Mapgap}

The following simple lemma relates the dynamics of $T_\be$ to
the original problem. Let
$$
\De_\be=[1/\be,1/\be(\be-1)].
$$

\begin{lemma}\cite{GS2}
\[
\A'_\be=\left\{x\in\left(0,\frac1{\be-1}\right) :
T_\be^n(x)\notin\De_\be\,\textrm{for all}\,\, n\ge0\right\}.
\]
\end{lemma}
\begin{proof} Note first that if $x\in[0,1/\be)$, then necessarily
$\e_1$ in (\ref{beta-exp}) is 0 and if
$x\in(1/\be(\be-1),1/(\be-1))$, then it is necessarily 1. If it
belongs to $\De_\be$, then there is always at least two
different representations.

Let $S_\be$ denote the shift on $\U_\be$ and
$R_\be=T_\be|_{\A_\be}$. Then the following commutative diagram
takes place:
\[
\begin{CD}
\U_\be @>{S_\be}>> \U_\be \\
@V{\pi_\be}VV @VV{\pi_\be}V \\
\A_\be @>{R_\be}>> \A_\be
\end{CD}
\]
Thus, $T_\be$ acts as a shift on sequences providing unique
representations. Indeed, it is either $x\mapsto\be x$ or $\be
x-1$ depending on whether we have $\e_1=0$ or 1, which is
precisely the shift in (\ref{beta-exp}). Hence for a
$T_\be$-orbit to stay out of $\De_\be$ at any iteration is the
same as keeping the representation (\ref{beta-exp}) unique.
\end{proof}

The important problem now is to describe the topological and
ergodic properties of the shift $S_\be$. The following theorem
summarizes all we know at present about them.

\begin{thm} \cite{GS2}
\begin{enumerate}
\item For every $\be\in(\be_*,2)$ the subshift $S_\be$ is {\bf
essentially transitive}, i.e., has a unique transitive component
of maximal entropy.
\item The shift $S_\be$ is a subshift of finite type for a.e.
$\be$.
\item For a.e $\be\in(\be_*,2)$ the subshift $S_\be$ is
intrinsically ergodic and metrically isomorphic to a {\bf
transitive} SFT.
\item The function $\be\mapsto h_{top}(\si_\be)$ is continuous
(but not H\"older continuous) and constant a.e. Every interval
of constancy is naturally parameterized by an algebraic integer
of a certain class.
\end{enumerate}\label{erg-prop}
\end{thm}

The open question is whether the shift $S_\be$ is ``as good as"
the beta-shift $\si_\be$ (see above). Namely, we {\em
conjecture} that for {\bf any} $\be\in(\be_*,2)$ it is
intrinsically ergodic and its natural extension is Bernoulli.

\begin{rmk} The above family of maps with gaps as a dynamical
object might look artificial: we change not only the slope but
the gap as well. However, if one alters the size of gap only
(which is more conventional), then the result on the symbolic
level will be essentially the same. For example, let
$Tx=2x\bmod1$ and $D_\de=[\de,1-\de]$ for $\de\in(0,1/2)$. Let
now
$$
\mathcal K_\de=\{x\in(0,1):T^n(x)\notin D_\de,\ n\ge0\}.
$$
It is known from the physical literature that $\dim_H(\mathcal
K_\de)>0$ if and only if $\de>\sum_1^\infty
\mathfrak{m}_n2^{-n}=0.412\dots$ (this has been apparently shown
independently in \cite{BorMc, ZB}). From the above results this
claim follows almost immediately. Sketch of the proof of this
fact is as follows: let $\mathbf a$ denote the binary expansion
of $2\de$; then in terms of the full 2-shift the set $\mathcal
K_\de$ is defined by (\ref{Uq}); the only difference with
$\U_\be$ is that the sequence $\mathbf a$ does not necessarily
satisfy the condition from Theorem~\ref{Pabe}~(\ref{P1}), but
this is easy to deal with. So, the shifted Thue-Morse sequence
is critical as well, which leads to the result in question. For
a more general case see \cite{GS2}.
\end{rmk}

\subsubsection{$\be>2$} Assume now that $\be\in(N,N+1)$
for some $N\ge2$; we have similar results with some natural
analog of the Thue-Morse sequence. Namely, let $\rho$ denote the
following substitution (morphism):
\[
a\to ac,\,\ b\to ad,\,\, c\to da,\,\, d\to db.
\]
Then
\[
\rho^\infty(d)=:(\mathfrak{w}_n)_1^\infty=dbabacdb\dots
\]
The sequence $(\mathfrak{w}_n)_1^\infty$ is related to the
Thue-Morse sequence in the following way:
\[ \si(\mathfrak{m})=
\underbrace{1101}_{d}\,\,\underbrace{0011}_{b}\,\,
\underbrace{0010}_{a}\,\,\underbrace{1101}_{d}\,\,
\underbrace{0010}_{a}\,\,\underbrace{1100}_{c}\,\,
\underbrace{1101}_{d}\,\,\underbrace{0011}_{b}\dots
\]
Let $P_n:\{0,1\}\to\{n-1,n\}$ be defined by
\[
P_n(0)=n-1,\ P_n(1)=n
\]
and $Q_n:\{a,b,c,d\}\to\{n-1,n,n+1\}$:
\[
Q_n(a)=n-1,\ Q_n(b)=Q_n(c)=n,\ Q_n(d)=n+1.
\]
Then the critical value $x=\be_*^{(N)}$ analogous to $\be_*$ is
given by the equations
\[
1=\sum_{k=1}^\infty P_n(\m_{k+1})x^{-k},\quad N=2n,\ n\ge1
\]
or
\[
1=\sum_{k=1}^\infty Q_n(\mathfrak w_k)x^{-k},\quad N=2n+1,\
n\ge1.
\]
More precisely, let $\A_\be^{(N)}$ denote the set of $x$ which
have a unique $\be$-expansion with the digits $0,1,\dots,N-1$.

\begin{thm} \cite{GS3}
\begin{enumerate}
\item The number $\be_*^{(N)}$ is the smallest $\be$ for which
$x=1$ has a unique $\be$-expansion with the digits
$0,1,\dots,N-1$.
\item The set $\A_\be^{(N)}$ has positive Hausdorff dimension if
and only if $\be>\be_*^{(N)}$. If $\be<\be_*^{(N)}$, then it is
at most countable.
\end{enumerate}
\end{thm}

\begin{rmk} The method of proving this type of theorems comes
from low-dimensional dynamics: it is called {\em
renormalization}. For more detailed results in this direction
and discussion see \cite{GS2, GS3}.
\end{rmk}

\subsection{Intermediate beta-expansions}\label{intermediate}
Assume for simplicity that $\be<2$ and have another look at the
map with a gap $T_\be$ defined above. More precisely, let
\begin{equation}
T'_\be(x)=\begin{cases}
\be x,& x\in\bigl[0,\frac1{\be(\be-1)}\bigr]\\
\be x-1,& x\in \bigl[\frac1{\be},\frac1{\be-1}\bigr].
          \end{cases}\label{multiv}
\end{equation}

Thus, we have a multivalued map on the middle interval
$\De_\be$, and in order to get a ``normal" (single-valued) map,
one needs to make a choice for every $x\in\De_\be$. From this
point of view the $\be$-shift corresponds to the choice of the
lower branch and the ``lazy" $\be$-shift -- the higher one for
every $x$. On the other hand, as we have seen, if one removes
the middle interval completely, then this corresponds to the
unique expansions and eventually to the map that acts on at most
a Cantor set, namely, $T_\be|_{\A_\be}$.

In \cite{DK} K.~Dajani and C.~Kraaikamp suggested the idea of
considering ``intermediate" beta-expansions, or how they called
them, $(\be,\al)$-expansions. More precisely, let $\al$ be our
parameter, $\al\in[0,(2-\be)/(\be-1)]$. We choose the upper
branch if the ordinate is less than $1+\al$ and the lower branch
otherwise, then restrict the resulting map to the interval
$[\al,1+\al]$. Thus, we get
\[
T_{\be,\al}(x)=\begin{cases} \be x,&
x\in\bigl[\al,\frac{1+\al}{\be}\bigr) \\
\be x-1,& x\in \bigl[\frac{1+\al}{\be},1+\al\bigr].
               \end{cases}
\]
(see Figure~\ref{beta-alpha}).

\pic{9truecm}{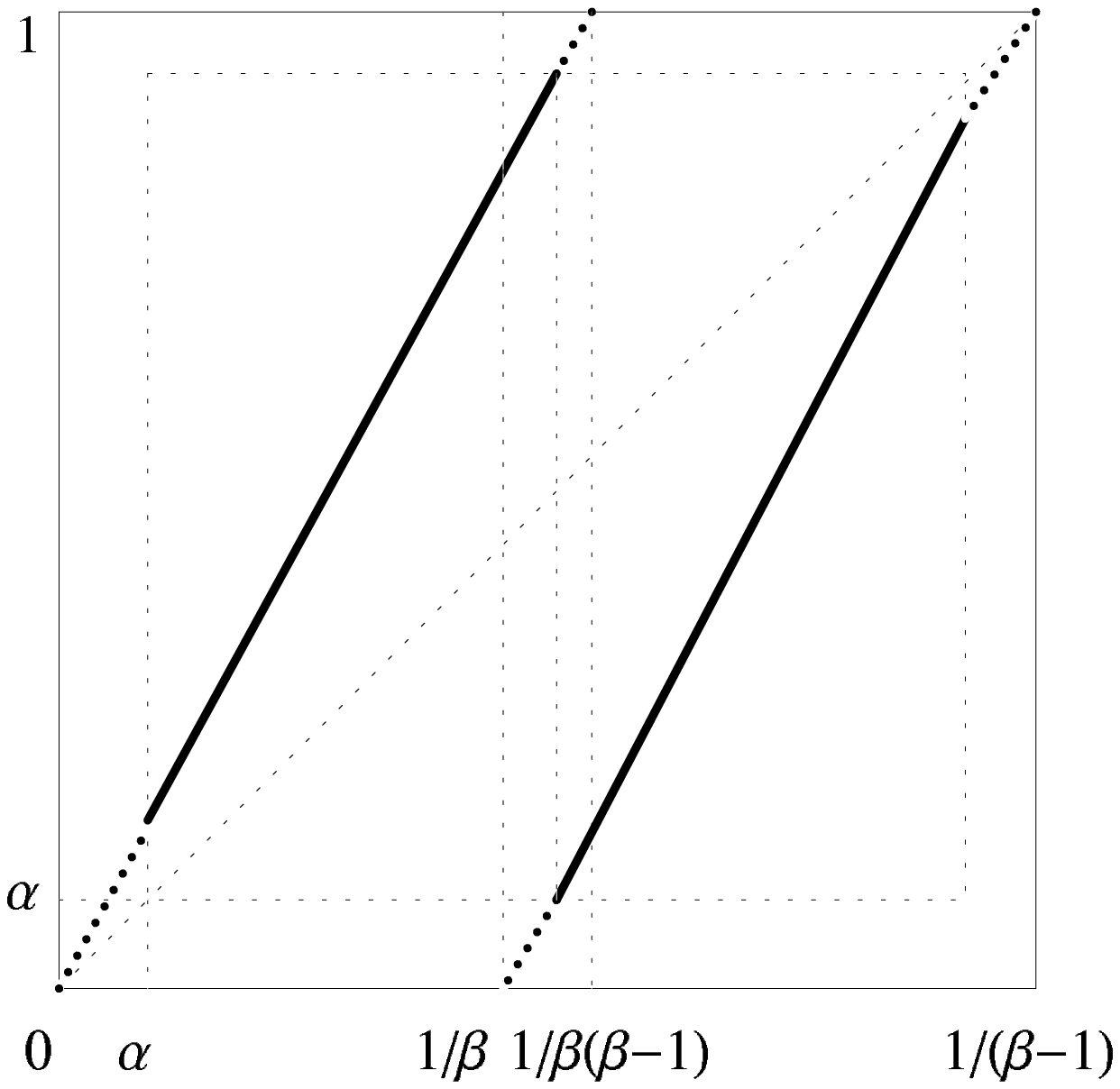}{The pattern for the
$(\be,\alpha)$-expansion}{beta-alpha}

As is easy to see, $T_{\be,\al}$ is isomorphic to the map
$S_{\be,\ga}:[0,1)\to[0,1)$ acting by the formula
\[
S_{\be,\ga}(x)=\be x+\ga\bmod1,\quad \ga=(\be-1)\al\in[0,2-\be].
\]

The ergodic properties of the family $S_{\be,\ga}$ are well
studied, see, {\em e.g.}, \cite{Pa2, Hof0, FlLag}. In
particular, $S_{\be,\ga}$ is ergodic with respect to the
Lebesgue measure, which is equivalent to the measure of maximal
entropy.

\subsection{The realm of beta-expansions} In the previous
subsections we have studied the beta-expansions from the
viewpoint of choosing a specific representation for a given $x$
in (\ref{beta-exp}) (and for $x\in\A'_\be$ this choice was
unique). We believe it is interesting to study the space of {\bf
all} possible representations of a given $x$, which this
subsection will be devoted to.

So, let $\be>1$ be fixed, and
$\Si_q=\prod_1^\infty\{0,1,\dots,q-1\}$ for some fixed $q\ge2$.
We define for any $x\ge0$,
\[
\R_{\be,q}(x):=\left\{\e\in\Si_q :
x=\sum_{n=1}^\infty\e_n\be^{-n}\right\}.
\]

In the case when $\be$ is non-integer, and $q=[\be]+1$, we will
simply write $\R_\be(x)$.

Not much is known about the general case yet. The following
result was obtained by the author \cite{generic} by elementary
means.\footnote{As with \cite{GS1}, this is a part of the
author's {\em credo}: an elementary result deserves an
elementary proof.} In this paper we would like to present an
ergodic-theoretic proof that may possibly start the whole new
line of research in this area (see Remark~\ref{branch} below).
The author is grateful to V.~Komornik for his help with the
history of the issue.

\begin{thm} \label{gencase} For any $\be\in(1,2)$ the cardinality
of $\R_{\be}(x)$ is the continuum for a.e. $x\in(0,1/(\be-1))$.
\end{thm}
\begin{proof} Recall that if $\be<G$, then {\bf every} point $x$
is known to have $2^{\aleph_0}$ $\be$-expansions
\cite[Theorem~3]{ErdJooKom}. So, let $\be\ge G$.

By the above, it suffices to show that
\begin{equation}\label{unionn}
\textrm{card}\,\left(\bigcup_{n=1}^\infty
(T_\be')^n(x)\right)=2^{\aleph_0}
\end{equation}
for a.e. $x\in(0,1/(\be-1))$ (where $T_\be'$ is the multivalued
map given by (\ref{multiv})). Let
\[
T_0(x)=
\begin{cases}
\be x,& x\in\bigl[0,\frac1{\be(\be-1)}\bigr)\\
\be x-1,& x\in\bigl[\frac1{\be(\be-1)},\frac1{\be-1}\bigr)
\end{cases}
\]
and
\[
T_1(x)=
\begin{cases}
\be x,& x\in\bigl[0,\frac1{\be}\bigr)\\
\be x-1,& x\in\bigl[\frac1{\be},\frac1{\be-1}\bigr)
\end{cases}
\]
(we omit the index $\be$ to simplify our notation). Then
$T_\be'x=T_0x\cup T_1x$ for any $x\in\De_\be$, whence any point
$y$ of the union in (\ref{unionn}) is of the form
\begin{equation} x_n=T_0^{k_n}T_1^{k_{n-1}}\dots
T_1^{k_2}T_0^{k_1}(x), \label{xnn}
\end{equation}
where $k_j\ge0, 1\le j\le n$. We are going to show that for a
Lebesgue-generic $x$ and every $x_n$ of the form (\ref{xnn})
there exists $k_{n+1}\in\BN$ such that
$T_\be^{k_{n+1}}(x_n)\in\De_\be$ (where $T_\be$ is given by
(\ref{Tbetaa})), i.e., the branching for the multivalued map
$T_\be'$ has the form of the binary tree for a.e. $x$.

Fix the vector $(k_1,\dots,k_n)\in\BZ_+^n$ and let within this
proof
\[
T:=T_0^{k_n}T_1^{k_{n-1}}\dots T_1^{k_2}T_0^{k_1}.
\]
Following the canonical proof of the corollary to the Poincar\'e
Recurrence Theorem (saying that a generic point returns to a
given set of positive measure infinitely many times), we will
show that for a.e. $x$ there exists $m\in\BN$ such that
$T^m(x)\in\De_\be$. This will prove the claim of our Theorem:
let $E(k_1,\dots,k_n)$ denote the generic set in question; then
\[
E:=\bigcap_{\substack{n\ge1\\(k_1,\dots,k_n)\in\BZ_+^n}}
E(k_1,\dots,k_n)
\]
will be a sought set of full measure.

So, it suffices to show that
\[
\mathcal L\left(\bigcup_{j=1}^\infty T^{-j}\De_\be\right)=1
\]
(where $\mathcal L$ is the {\bf normalized} Lebesgue measure on
$[0,1/(\be-1)]$). But since $T_0$ and $T_1$ are both extensions
of $T_\be$, we have $T^{-1}A\supset T_\be^{-p}A$ for every
measurable $A$ with $p=k_1+\dots+k_n$. It is thus left to show
that
\begin{equation}
\mathcal L\left(\bigcup_{j=1}^\infty T_\be^{-pj}\De_\be\right)=1
\label{cupp}
\end{equation}
for every $p\ge1$. Let
\[
\wt T_\be(x)=
\begin{cases}
\be x,& x\in\bigl[0,\frac1{\be}\bigr]\\
1,& x\in\bigl(\frac1{\be},\frac1{\be(\be-1)}\bigr)\\
\be x-1,& x\in\bigl[\frac1{\be(\be-1)},\frac1{\be-1}\bigr].
\end{cases}
\]
Thus, $\wt T_\be$ is a ``real map"; however, it is obvious that
preimage-wise the maps $T_\be$ and $\wt T_\be$ are the same.
Moreover, the normalized Lebesgue measure on $[0,1/(\be-1)]$ is
quasi-invariant under $\wt T_\be$; similarly to \cite{Re} and
\cite{Pa2}, one can easily show that there exists a unique $\wt
T_\be$-invariant measure $\mathcal L_\be$ which is equivalent to
the Lebesgue measure. Thus, the measure $\mathcal L_\be$ is
$(\wt T_\be)^p$-invariant as well, which by the Poincar\'e
Recurrence Theorem implies (\ref{cupp}), and we are done.
\end{proof}
\begin{rmk} It is easy to generalize this theorem to the case of
an arbitrary non-integer $\be>1$ and arbitrary $q>[\be]$. We
leave the details for the reader as a simple exercise (note: it
is sufficient to consider the case $q=[\be]+1$).
\end{rmk}
\begin{rmk}\label{branch} The branching for $T_\be'$ described
in the proof of the theorem, may be apparently studied in a more
quantitative way. In particular, we {\em conjecture} that not
only the cardinality of $\R_\be(x)$ is the continuum for a
generic $x$ but its Hausdorff dimension in the space $\Si$
(provided with the natural ($=$ binary) metric) is positive.
This might be possibly shown by studying the average return
times to the set $\De_\be$ for the multivalued map $T_\be'$ and
will be considered elsewhere.
\end{rmk}

\begin{exam}\label{ex12} Let $\be=G$; it is shown in
\cite[Appendix~A]{SV1} that $\R_G(x)$ is always a continuum
unless $x=nG\bmod1$ for some $n\in\BZ$. For example,
$\R_G(1/2)=\prod_1^\infty\{011,100\}$, whence it is indeed a
continuum, and its Hausdorff dimension equals $\frac13$.
\end{exam}

An important special case is when $\be$ is a Pisot number (see
Section~\ref{greedy} for definition). Then one may study the
combinatorics of the infinite space by means of the
combinatorics of its finite approximations. More precisely, let
us call two sequences $\e$ and $\e'$ from $\Si_q$ {\em
equivalent} if
$\sum_1^\infty\e_k\be^{-k}=\sum_1^\infty\e'_k\be^{-k}$. There
are a number of results about counting the cardinality of
equivalence classes; most of them deal with random matrix
products \cite{Lal1, Lal2, OST}. In particular, it is shown in
\cite{Lal2} that for a generic 0-1 sequence $(x_n)_1^\infty$ the
cardinality of the equivalence class of $(x_1\dots x_n)$ grows
with $n$ exponentially with an exponent greater than 1. The
value of this exponent is given by the upper Lyapunov exponent
of the random matrix product in question.

One special case is however worth mentioning on its own. Let
$\be=G$ (the golden ratio) and $q=2$. We call a finite 0-1 word
a {\em block} if it begins with 1 and ends by an even number of
0's. As is easy to show by induction, every block has the
following form:
\[
B=1(01)^{a_1}(00)^{a_2}\dots (00)^{a_r}\quad \mathrm{or}\quad
B=1(00)^{a_1}(01)^{a_2}\dots (00)^{a_r}.
\]
Hence each block $B$ is parameterized in a unique way by the
natural numbers $a_1,\dots, a_r$. We will write
$B=B(a_1,\dots,a_r)$. Let $\c(w)$ denote the cardinality of the
set of all 0-1 words equivalent to $w$.

In 1998 the author together with A.~Vershik proved the following

\begin{thm} \cite{SV1}
\begin{enumerate}
\item Let $w=B_1\dots B_k$, where $B_j$ is a block for all $j$.
Then the space of all 0-1 words equivalent to $w$ splits into
the Cartesian product of the equivalence classes for $B_j$ for
$j=1$ to $k$, and therefore, the function $\c$ is blockwise
multiplicative:
\[
\c(B_1\dots B_k)=\c(B_1)\dots\c(B_k).
\]
\item The cardinality of a block is given by the formula
\[
\c(B)=p_r+q_r,
\]
where $\frac{p_r}{q_r}=[a_1,\dots,a_r]$ is the continued
fraction.
\end{enumerate}\label{block}
\end{thm}

The proof given in \cite{SV1} is based on an induction argument;
however, now it is clear that there exists a more direct and
elegant way of proving this result. Let $X=X_G$ which in this
case is the Markov compactum of all 0-1 sequences without two
consecutive 1's and let $g:\prod_1^\infty\{a,b,c\}\to X$ act by
the following rule:
\[
g(a)=00,\ g(b)=010,\ g(c)=10
\]
and then by concatenation.
\begin{lemma} The map $g$ is a bijection.
\end{lemma}
\begin{proof} It is a straightforward check that $g^{-1}$ is well
defined: any sequence from $X$ can be split in a unique way into
the concatenation of the blocks $00, 010$ and $10$.
\end{proof}

Let $w$ be a finite word in $X$ and $g^{-1}w=j_1\dots j_m$. In
\cite{OST} it was shown that
\[
\c(w)=\begin{pmatrix} 1&0\end{pmatrix}P_{j_1}P_{j_2}\dots
P_{j_m}\begin{pmatrix} 1\\0\end{pmatrix},
\]
where
\[
P_a=\begin{pmatrix} 1 & 1 \\ 0 & 1
\end{pmatrix},\quad P_b=\frac12\begin{pmatrix} 1 & 1 \\ 1 & 1
\end{pmatrix},\quad P_c=\begin{pmatrix} 1 & 0 \\ 1 & 1
\end{pmatrix}.
\]
This immediately yields both parts of Theorem~\ref{block}, in
view of the well-known relation between the matrix products of
$P_a$ and $P_c$ and the finite continued fractions. The details
are left to the reader.

The main consequence of Theorem~\ref{block} is the existence of
the map called the {\em goldenshift} which acts on sequences
from $X$ starting with 1 as the shift by the length of the first
block. This goldenshift is well defined a.e. and has a number of
important properties that reveal a lot of information about the
Bernoulli convolution parameterized by the golden ratio. For
details see \cite[\S\S2, 3]{SV1}.

\begin{rmk} In \cite[Appendix~A]{SV1} A.~Vershik and the
author completely described all the possible patterns for
$\R_G(x)$ if $x\in(0,1)$. Loosely speaking, these possibilities
are as follows: either $\R_G(x)$ is a Cartesian product (this is
Lebesgue-generic -- see also Example~\ref{ex12}) or
$\R_G(x)=\Xal'$ for a certain $\al=\al(x)$, where $\Xal'$ is the
``rotational" Markov compactum described in the next section.
For more general Pisot numbers the analog of this theorem seems
to be a delicate and interesting problem.
\end{rmk}

\begin{rmk} The combinatorics of the ``integral" case
$\be\in\BN$ is studied in detail by J.-M.~Dumont, A.~Thomas and
the author in \cite[\S6]{DST}. In that case the cardinality
function can be represented in terms of random matrix products
as well. Specifically, if $\be=2$ and $q=3$, then these matrices
are precisely $P_a$ and $P_c$.
\end{rmk}

\begin{rmk} There exists a class of singular measures (Bernoulli
convolutions) that are based on the combinatorics we have just
described. For more details see, {\em e.g.}, \cite{AlZa, SV1}
for the case of $\be=G$ and \cite{E, Ga, Lal1, Lal2, OST} for
more general cases of Pisot numbers (see also the survey article
\cite{PSS} for a general overview).
\end{rmk}

\section{Rotational expansions}\label{RE}

In this section we are going to describe the model appeared
first in \cite{Ver82} as an application of the general theorem
by A.~Vershik on adic realization -- see
Theorem~\ref{metric_iso} in Appendix, and in a more arithmetic
form -- in the joint paper \cite{VS} by A.~Vershik and the
author. Also, we will present another model which deals with
more conventional base for arithmetic expansions.

\subsection{General constructions}\label{RE-1}
\subsubsection{First model}
The problem we are going to consider in this section, is
arithmetic codings of an irrational rotation of the circle
$\BR/\BZ$ which we will identify with the interval $[0,1)$. Let
$\al\in(0,1/2)\setminus\BQ$ be the angle of rotation (if it is
greater than $1/2$, simply take $1-\al$), and let
$\Ral(x)=x+\al\bmod1$. Let the regular continued fraction
expansion of $\al$ be
\[
\al=\cfrac{1}{a_1+\cfrac{1}{a_2+\cfrac{1}{a_3+\dots}}}
\]
and $(p_n/q_n)_1^\infty$ be the sequence of convergents with
$p_1=0,p_1=1,q_1=1,q_2=a_1$ and
$p_{n+1}=a_np_n+p_{n-1},q_{n+1}=a_nq_n+q_{n-1}$. Since
$\al<1/2$, we have $a_1\ge2$. Put $r_1=a_1,r_n=a_n+1, n\ge2$
and,
\[
M^{(n)}=M^{(n)}(\al):=\begin{pmatrix} 1 & \dots& 1&1\\
1&\dots& 1&0\\ \vdots& \vdots& \vdots&\vdots \\
1&\dots &1&0\end{pmatrix},
\]
where the number of rows in $M^{(n)}$ is $r_n$, and the number of
columns is $r_{n+1}$. Let $D_n=\{0,1,\dots,r_n-1\}$ be endowed
with the natural ordering, i.e., $0\prec1\prec\cdots\prec r_n-1$.
Let now $\Xal$ denote the Markov compactum determined by the
sequence of matrices $(M^{(n)}(\al))_1^\infty$ and $\T$ denote the
adic transformation on $\Xal$ with respect to the natural ordering
(see Appendix for the definitions).

Note first that $\T$ is well defined everywhere with the exception
of two ``maximal" sequences: $(a_1-1,0,a_3,0,a_5,\dots)$ and
$(0,a_2,0,a_4,0,\dots)$. Thus, if we exclude the {\em cofinite}
sequences (i.e., those whose tail is $(a_n,0,a_{n+2},0,\dots)$ for
some $n$), the positive part of the orbit of $\T$ will be always
well defined. To enable the whole trajectory of $\T$ to be well
defined, one has to remove the finite sequences as well.

We claim that the map $\T$ is metrically isomorphic to $\Ral$
and are going to present the conjugating map. Let
$\psi_\al:\Xal\to[0,1)$ be defined as
\begin{equation}
\psi_\al(x_1,x_2,\dots):=\al+\sum_{n=1}^\infty
x_n(-1)^{n+1}\al_n,
\label{conjmap}
\end{equation}
where $\al_n=|q_n\al-p_n|=\|q_n\al\|$ (here $\|\cdot\|$ stands
for the distance to the nearest integer).

\begin{thm} \cite{VS}
\begin{enumerate}
\item The invertible dynamical system $(\Xal,\T)$ is uniquely
ergodic. The unique invariant measure $\nu_\al$ is Markov on
$\Xal$.
\item $\psi_\al$ is continuous and one-to-one except the
cofinite sequences.\label{Ss}
\item The map $\psi_\al$ metrically conjugates the
automorphisms $(\Xal,\nu_\al,\T)$ and $([0,1),\Leb,\Ral)$, where
$\Leb$ stands for the Lebesgue measure on the unit interval.
\end{enumerate}\label{VS1}
\end{thm}

\begin{rmk} As is well known, $\al_n$ decays at a very fast rate
as $n$ goes to the infinity, namely, $\al_n=O(q_{n+1}^{-1})$
(see, {\em e.g.}, \cite{Khin}). Thus, the $n$'th term of the sum
in (\ref{conjmap}) is $O(q_n^{-1})$, i.e., decays at worst
exponentially.
\end{rmk}

\begin{rmk} Expansion~(\ref{conjmap}) was considered for the
first time by Y.~Dupain and V.~Sos \cite{Sos}; they also proved
Theorem~\ref{VS1}~(\ref{Ss}), the fact A.~Vershik and the author
were unaware of when writing \cite{VS}. This however hardly
undermines Theorem~\ref{VS1} as it appeared in \cite{VS}, because
the (most important) dynamical meaning of the rotational expansion
was new.
\end{rmk}

\begin{rmk}
A simple way to obtain these expansions is as follows: let
$N\in\BN$ and $(q_n)_1^\infty$ serve as a ``base" for
representations in the sense of \cite{Fra}. Then there is a
unique representation of $N$ in the form
\begin{equation}
N=1+\sum_k x_kq_k,\label{N=}
\end{equation}
where $(x_1,x_2,\dots)$ is a finite sequence in $\Xal$. All one
has to do to get (\ref{conjmap}) is to make a {\em profinite
completion} of (\ref{N=}) using the fact that the sequence
$(N\al\bmod1)_1^\infty$ is dense in $[0,1)$ and the formula
$q_k-p_k\al=(-1)^{k+1}\al_k$. For more details see \cite[\S2]{VS}.
\end{rmk}

\begin{rmk} If one removes both the finite and cofinite
sequences from $\Xal$ and the $\T$-trajectory of 0 from $[0,1)$,
then $\psi_\al$ becomes a homeomorphism and thus, acts as a
conjugacy in the topological sense as well.
\end{rmk}

\subsubsection{Second model} The price we pay for the natural
ordering in the first model is that the base of the expansions
is not always positive. The second model we are going to
describe below, overcomes this problem but here there is the
price to pay as well: the ordering is rather unusual.
Apparently, it is impossible to take care of both issues
simultaneously -- this symbolizes the well-known fact that the
convergents $(p_n/q_n)<\al$ if $n$ is even and $>\al$ if $n$ is
odd.

Let $\al\in(0,1)\setminus\BQ$ and
\[
\mathcal{M}^{(n)}=\mathcal{M}^{(n)}(\al):=
\begin{pmatrix} 1 & 1&\dots& 1&1\\
\dots&\dots&\dots&\dots&\dots\\
1 & 1&\dots& 1&1\\
1&0&\dots& 0&0\\
\end{pmatrix},
\]
where the size of $\mathcal{M}^{(n)}$ is $(a_n+1)\times
(a_{n+1}+1)$. Let $\mathcal{D}_n=\{0,1,\dots,a_n\}$. We define
the {\em alternating} ordering on $\mathcal D_n$ as follows:
$0\prec1\prec\cdots\prec a_n$ if $n$ is odd and
$0\succ1\succ\cdots\succ a_n$ if $n$ is even.

Let now $\mathcal X'_\al$ denote the corresponding Markov
compactum, $\mathcal T'_\al$ -- the adic transformation on it
and $\psi'_\al:\mathcal X'_\al\to[0,1)$ be defined as follows:
\begin{equation}
\psi'_\al(x_1,x_2,\dots):=\sum_{n=1}^\infty x_n\al_n,
\label{conjmap2}
\end{equation}
where $\al_n$ are as above. We have the following analog of
Theorem~\ref{VS1}:

\begin{thm} \cite{VS}
\begin{enumerate}
\item The map $\mathcal T'_\al$ is well defined everywhere
except the sequence $(a_1,0,a_3,\linebreak[1] 0,\dots)$. Its
inverse is not well defined only at $(0,a_2,0,a_4,0,\dots)$.
\item The dynamical system $(\mathcal X'_\al,\mathcal T'_\al)$ is
uniquely ergodic. The unique invariant measure $\mu_\al$ is
Markov on $\Xal'$.
\item $\psi'_\al$ is continuous and one-to-one except the
cofinite sequences.
\item The map $\psi'_\al$ metrically conjugates the automorphisms
$(\Xal',\mu_\al,\T')$ and $([0,1),\Leb,\Ral)$.
\end{enumerate}
\end{thm}

\begin{rmk} The expansion (\ref{conjmap2}) is a special case of
the general class of {\em Cantor-Waterman expansions}
\cite{Wat}. More general systems of numeration are considered in
\cite{VS} as well but without clear dynamical meaning.
\end{rmk}

\begin{rmk}
An analog of (\ref{N=}) for the integers is as follows:
\begin{equation}
N=\sum_n x_n(-1)^n q_n,\label{N==}
\end{equation}
where $N\in\BZ$ (not necessarily nonnegative!) and $(x_n)$ is a
finite sequence from $\Xal'$ \cite{VS}. Thus, the alternating base
is not for the reals but for the integers in the second model.
Comparing (\ref{conjmap}) with (\ref{conjmap2}) and (\ref{N=})
with (\ref{N==}), we see that the two models are in a way dual.
\end{rmk}

\subsection{Probabilistic properties of the ``digits"} In this
subsection we are going to mention briefly the results from
\cite{S95} on Laws of Large Numbers (LLN and SLLN) and the Central
Limit Theorem (CLT) for the sequences of digits for the expansions
considered above. We will consider the metric space
$(\Xal',\mu_\al)$; the results for $(\Xal,\nu_\al)$ are very
similar, so we will omit them.

First, the initial distribution for $\mu_\al$ is as follows:
\[
\mu_\al(x_1=i_1)=
\begin{cases}
\al,& x_1<a_1 \\ \al_2,& x_1=a_1,
\end{cases}
\]
the transition probabilities are given by
\begin{equation}
\mu_\al\big(x_n=i_n\mid x_{n-1}=i_{n-1}\big)=
\begin{cases}
\frac{\al_n}{\al_{n-1}},& i_{n-1}<a_{n-1},\quad i_n<a_n\\
\frac{\al_{n+1}}{\al_{n-1}},& i_{n-1}<a_{n-1},\quad i_n=a_n\\
1,& i_{n-1}=a_{n-1},\quad i_n=0\\
0,& \text{otherwise}.
\end{cases}
\label{trprob}
\end{equation}
Finally, the one-dimensional distributions are as follows:
\begin{equation}
\mu_\al\big(x_n=i_n)=
\begin{cases}
\big(q_{n-1}+q_n\big)\al_n, & i_n=0\\
q_n\al_n, & 0<i_n<a_n\\
q_n\al_{n+1}, & i_n=a_n.
\end{cases}
\label{oddistr}
\end{equation}

The following theorem shows that, roughly speaking, if the
partial quotients $a_n$ of $\al$ do no grow ``too fast", then
most of the probabilistic laws hold. More precisely, we have

\begin{thm}\cite{S95}
\begin{enumerate}
\item If
\[
\sum^n_{k=1}a^2_k=o(n^2),\quad n\to\infty,
\]
then the LLN holds for $(\Xal',\mu_\al)$.
\item The condition
$$
\sum^\infty_{n=1}\frac{a^2_n}{n^2}\ln^2 n<+\infty
$$
is sufficient for the validity of SLLN for $(\Xal',\mu_\al)$.
\item Finally, if the partial quotients for $\al$ are uniformly
bounded, then the CLT holds for $(\Xal',\mu_\al)$.
\end{enumerate}
\end{thm}

\begin{rmk} None of these conditions is Lebesgue-generic for
$\al$. We believe all of them can be improved but not
significantly.
\end{rmk}

\subsection{Unique rotational expansions}\label{URE} Following
the pattern of Section~\ref{un}, it is interesting to study the
combinatorics of expansions~(\ref{conjmap2}) with the lifted
Markov restrictions (again, we will not consider the first
model, where all the results is very similar and leave it to the
interested reader). The results presented below are original
(though some steps in this direction have been undertaken in
\cite{S95-2}).

Let $Z_\al=\prod_{n=1}^\infty\{0,1,\dots,a_n\}$ and
\[
\V:=\left\{x\in(0,1)\mid\exists !\ (x_1,x_2,\dots)\in Z_\al :
x=\sum_{n=1}^\infty x_n\al_n\right\}.
\]
Our goal will be to study the properties of $\V$.

Let us recall that $\al_{n-1}=a_n\al_n+\al_{n+1}$. Hence the
triples $x_{n-1}=1,x_n=0,x_{n+1}=0$ and
$x_{n-1}=0,x_n=a_n,x_{n+1}=1$ give the same value in
(\ref{conjmap2}) provided all the other digits are the same. In
a way, this claim is invertible, and this is what the proof will
be based upon.

More precisely, in \cite{S95-2} it is shown that if $x\in(0,1)$
has at least two different representations in the form
(\ref{conjmap2}), then in its {\em canonical} representation
(the one with the digits in $\Xal'$) there exists $n\in\BN$ and
a triple $(i_{n-1},i_n,i_{n+1})$ with $i_{n-1}>0, i_n=0$ and
$i_{n+1}<a_{n+1}$. We will call such triples {\em replaceable}.
The question is, whether replaceable triples are generic with
respect to the Lebesgue measure.\footnote{The condition looks
like something shift-invariant but there is no suitable ergodic
theorem here, of course - the compactum $\Xal'$ is
non-stationary! (unless $a_n\equiv a$ for all $n$)}

We denote
\[
\V'=(\psi'_\al)^{-1}(\V),
\]
i.e., the set of admissible sequences providing unique rotational
representations.

\begin{lemma} Let $(x_1,x_2,\dots)\in\V'$. Then
\begin{enumerate}
\item if $x_n=0$ for some $n$, then necessarily
$x_{n-1}=\dots=x_1=0$ as well;
\item it is impossible that $x_n=a_n$ for some $n$.
\end{enumerate}\label{tt}
\end{lemma}
\begin{proof} It suffices to prove (1), because $x_n=a_n$
implies $x_{n+1}=0$, which would contradict (1). Assume $x_n=0$;
if $x_{n-1}>0$, then the triple $(x_{n-1},x_n,x_{n+1})$ will be
replaceable unless $x_{n+1}=a_{n+1}$. But then again, we have
$x_{n+2}=0$, which leads to the same problem! Since the tail
$(x_n=a_n,0,a_{n+2},0,\dots)$ is not admissible, we are done.
\end{proof}

\begin{prop}\label{empty} $\V=\emptyset$ if and only if
\begin{equation}
\#\{n:a_n=1\}=+\infty.\label{an=1}
\end{equation}
\end{prop}
\begin{proof} (1) Assume $a_{n_k}=1$ for $k=1,2,\dots$ and
$(x_1,x_2,\dots)\in\V'$. Then for each $k$ we have a choice
between $x_{n_k}=0$ and $x_{n_k}=1$. The latter is impossible by
Lemma~\ref{tt}, while the former leads to $x_j\equiv0$ for all
$j\le n_k$, which in turn leads to $x_j\equiv0$ for all $j\in\BN$.
This is a contradiction, because $x>0$.\newline (2) If the number
of $a_n=1$ is finite, we set
\begin{equation}
n_0=\sup\,\{n:a_n=1\}.\label{n0}
\end{equation}
Then the sequence with $x_j=0$ for $1\le j\le n_0$ and $x_j=1$
otherwise, is a unique rotational representation.
\end{proof}
\begin{rmk} The condition~(\ref{an=1}) is Lebesgue-generic for
$\al$. So, for a typical $\al$ our set is empty.
\end{rmk}

The following result may be regarded as an analog of
Proposition~\ref{Leb=0} for the rotational expansions. The
crucial difference is that it is not true that for every
irrational $\al$ the set $\V$ has zero Lebesgue measure -- this
depends on how fast the partial quotients grow.

\begin{thm} The set $\V$ has Lebesgue measure zero if and
only if
\begin{equation}
\sum_{n=1}^\infty\frac1{a_n}=+\infty.\label{star}
\end{equation}
\end{thm}
\begin{proof} Assume first that (\ref{star}) is {\bf not} satisfied.
Then there exists only a finite number of $n$ such that $a_n=1$.
Let $n_0$ be given by (\ref{n0}) and
\[
K_\al:=\{(x_1,x_2,\dots)\in\Xal':0<x_n<a_n,\ n\ge n_0+1\}.
\]
By the above, each sequence from $K_\al$ is a unique
representation, whence it would suffice to show that
$\mu_\al(K_\al)>0$. This measure can be computed explicitly: by
(\ref{oddistr}) and (\ref{trprob}) and in view of $x_1=\dots=
x_{n_0}=0$ (see Lemma~\ref{tt}),
\[
\mu_\al(K_\al)=\al_{n_0}\cdot\prod_{n=n_0+1}^\infty
(a_n-1)\frac{\al_n}{\al_{n-1}}>
\al_{n_0}\prod_{n=n_0+1}^\infty\frac{a_n-1}{a_n+1}>0
\]
(the first inequality follows from the fact that
$(a_n+1)\al_n>\al_{n-1}$ and the second one is a consequence of
failing of (\ref{star})).

Assume now that (\ref{star}) holds. Our goal is to show that
$\mu_\al(\V')=0$, and our first remark consists in the
observation that by Proposition~\ref{empty}, it suffices to
consider $\al$ such that $n_0<\infty$. Furthermore, each
sequence from $\V'$ that contains at least one zero, is of the
form $(0,0,\dots,0,x_n,x_{n+1},\dots)$, where $0<x_j<a_j$ for
$j\ge n$ (see Lemma~\ref{tt}). Thus, we come again to the set
similar to $K_\al$ -- see above. We have
\[
\mu_\al(\V')= \sum_{k=0}^\infty\al_k\prod_{n=k+1}^\infty
(a_n-1)\frac{\al_n}{\al_{n-1}}\le\sum_{k=0}^\infty\al_k
\prod_{n=k+1}^\infty\frac{a_n-1}{a_n}=0,
\]
because by (\ref{star}), each infinite product in the
last-mentioned sum equals 0.
\end{proof}
\begin{rmk} In \cite{S95-2} it is shown that if (\ref{star})
is not satisfied, then the image of the uniform measure on
$Z_\al$ (i.e., $\prod_{n=1}^\infty\{1/a_n,\dots,1/a_n\}$) under
the map $\psi_\al'$ given by (\ref{conjmap2}), is an absolutely
continuous measure. In the opposite direction the result is
incomplete: apart from (\ref{star}) for this measure to be
singular, there is one (apparently, parasite) condition, which
at the time we have not been able to get rid of. Of course, if,
for instance, $a_n\equiv a$ for all $n\ge1$, then the measure in
question is singular, which is the famous Erd\"os Theorem
\cite{E}.
\end{rmk}

What is left if we wish to follow the pattern of
Section~\ref{un}, is the Hausdorff dimension of $\V$ when
(\ref{star}) is satisfied.

\begin{prop} Assume that the number of $n$ such that $a_n=1$, is
finite. Then the cardinality of $\V$ is the continuum if and
only if the tail of $(a_1,a_2,\dots)$ is different from
$(2,2,2,\dots)$. Otherwise $\V$ is a finite set.\label{contin}
\end{prop}
\begin{proof} Let again $n_0$ be given by (\ref{n0}). If
$a_n\equiv2$ for $n\ge n_1$, then we must have $x_n\equiv1$ for
$n\ge n_1$. If, on the contrary, there exists a subsequence
$(m_k)$ such that $a_{m_k}\ge3$, then we will have a choice of
$x_{m_k}=1$ or 2, which yields a continuum.
\end{proof}
\begin{rmk} Thus, here we also have some kind of monotonicity,
namely, the cardinality and Hausdorff dimension of $\V$ are
nondecreasing functions with respect to the partial quotients.
The difference with Section~\ref{un} is that $\V$ is never
infinite countable.
\end{rmk}

\begin{thm} Under the assumption of Proposition~\ref{contin},
the Hausdorff dimension of $\V$ is positive if and only if
\begin{equation}
\liminf_{n\to+\infty}\frac{\sum_{k=n_0+1}^n\log(a_k-1)}{\log
q_{n+1}}>0, \label{HD}
\end{equation}
where $n_0$ is given by (\ref{n0}).
\end{thm}
\begin{proof} Assume for the simplicity of notation that $a_n>1$
for all $n\ge1$. Let $\V^{(n)}$ denote the set of all cylinders of
length~$n$ in $\V'$. By the above, we will have the following
choice for $\V^{(n)}$: if $a_k=2$, then necessarily $x_k=1$;
otherwise $x_k\in\{1,2,\dots,a_k-1\}$. Hence by (\ref{oddistr})
and (\ref{trprob}),
\[
\mu_\al(\V^{(n)})=\al_n\prod_{k=1}^n(a_k-1)
\]
(all the transitional measures at the $k$'th step are the same),
and the condition for the positivity of the Hausdorff dimension
of $\V$ is a follows:
\[
\liminf_n\frac{\log\prod_{k=1}^n(a_k-1)}{-\log\al_n}>0,
\]
which, in view of the inequality $1/2<q_{n+1}\al_n<1$ is
equivalent to (\ref{HD}).
\end{proof}
\begin{cor} If $n_0$ given by (\ref{n0}) is less than infinity,
and
\[
\liminf_{n\to+\infty}\frac{\sum_{k=n_0+1}^n\log(a_k-1)}
{\sum_{k=n_0+1}^n\log(a_k+1)}>0,
\]
then $\dim_H(\V)>0$.
\end{cor}
\begin{proof} It suffices to use the relation
$q_{n+1}=a_nq_n+q_{n-1}$, from which it follows that
$q_{n+1}\le\prod_{k=1}^n(a_k+1)$.
\end{proof}

\begin{cor} If $a_n\le C$ for all $n\ge1$ and $n_0<\infty$,
then $\dim_H(\V)>0$ if and only if
\[
\De=\liminf_{n\to\infty}\frac1n\#\{1\le k\le n : a_k=2\}<1.
\]
\begin{proof} Again, for simplicity we assume that $a_n\neq1$ for
all $n\ge1$. We have: $\sum_{k=1}^n\log(a_k-1)\ge(1-\De)n\log2$,
and $\sum_{k=1}^n\log(a_k+1)\le n\log(1+C)$, whence by the
previous corollary, the ``if" part follows. The proof of the
``only if" part is left to the reader.
\end{proof}
\end{cor}

\section{Arithmetic codings of toral automorphisms}\label{AC}

This section is devoted to the arithmetic codings of hyperbolic
automorphisms of a torus. The idea of a coding is to expand the
points of a torus in power series in base its homoclinic point.
It was suggested by A.~Vershik in special cases \cite{Ver-ML,
Ver} and developed by the author and A.~Vershik in \cite{SV1,
SV2, VS2} in dimension~2 and in higher dimensions
(chronologically) -- by R.~Kenyon and A.~Vershik \cite{KenVer},
S.~Le~Borgne in his Ph.~D. Thesis \cite{Leb} and subsequent
works \cite{Leb1, Leb2}, K.~Schmidt \cite{Sch} and finally by
the author \cite{S1}.

\subsection{An important example: the Fibonacci automorphism} We
begin with the example that was studied in detail in 1991--92
and has eventually led to the theory described in the rest of
the section.

We are going to expose it just the way it appeared. The initial
motivation has come from the theory of $p$-adic numbers: let $p$
be a prime, and $Z_p$ denote the group of $p$-adic integers, i.e.,
one-sided formal series in powers of $p$:
\[
Z_p=\left\{\sum_{n=-\infty}^{-1} x_np^{-n} : 0\le x_n\le p-1,\
n\le -1 \right\}.
\]
Let $Q_p$ denote the field of $p$-adic numbers, i.e.,
\[
Q_p=\left\{\sum_{n=-\infty}^\infty x_np^{-n} \mid 0\le x_n\le
p-1,\ n\in\BZ,\ \exists N\in\BZ : x_n\equiv0, n\ge N\right\}.
\]
Thus, $Q_p$ is the space of two-sided $p$-adic expansions finite
to the right.\footnote{I have heard some people call them
``1.5-sided expansions". Informally, of course.} Finally, if one
considers the ``full-scale" two-sided $p$-adic expansions
\[
\mathcal S_p=\left\{\sum_{n=-\infty}^\infty x_np^{-n}\mid 0\le
x_n\le p-1,\ n\in\BZ\right\},
\]
then we obtain the $p$-adic solenoid.

\begin{question} What will all the above objects become if one
replaces $p$ by an algebraic unit $\be>1$ and the full $p$-adic
compactum -- by the two-sided $\be$-compactum $\wt X_\be$?
\end{question}

The obvious candidate to start investigation seemed
$\be=\frac12(1+\sqrt5)$, in which case, we recall, $\wt X:=\wt
X_\be$ is the set of two-sided 0-1 sequences without two
consecutive 1's (see Section~\ref{BE}). There is another good
reason for considering the golden ratio. Let $F_1=1,
F_2=2,\dots$ be the Fibonacci sequence; as was explained in
Section~\ref{RE}, every natural number $N$ has a unique
representation in base $(F_n)_1^\infty$ with the digits from
$X=X_\be$ -- see (\ref{N=}). Furthermore, as we know, the
profinite completion of (\ref{N=}) turns $\BN$ into $S^1$,
whence the analog of $Z_p$ is $S^1$. This suggests that unlike
the $p$-adic case, the {\em fibadic} case, as we will call it,
produces the topology of the real line instead of the $p$-adic
topology.

Recall also that the set of $p$-adic expansions as well as
$\be$-expansions finite to the {\bf left}, is simply $\BR_+$
(Lemma~\ref{half}). The situation with the spaces that involve
{\em formal} power series in base $\be$ (infinite to the left) is
completely different and strongly depends on $\be$, as we will see
below.

To deal with the problems regarding the formal power series, we
notice that in the $p$-adic case the key to the structure of
$Z_p, Q_p$ and $\mathcal S_p$ is just the following relation:
$pv_n=v_{n-1}$, where $v_n=p^{-n}$. In the fibadic case the
analog of this relation is
\begin{equation}
u_{n-1}=u_n+u_{n+1}.\label{fib}
\end{equation}
Thus, for instance, the analog of $Q_p$ is as follows:
\[
Q_\be:=\left\{\sum_{n=-\infty}^\infty \e_nu_n :
(\e_n)_{-\infty}^\infty\in\wt X,\ \e_n\equiv0,\ n\ge
N\,\,\mathrm{for}\,\,\mathrm{some}\,\, N\in\BZ \right\},
\]
where the sequence $(u_n)$ satisfies (\ref{fib}).

\begin{prop}\cite{Ver-ML} After identification of a countable
number of certain pairs of sequences $Q_\be$ becomes a field
isomorphic to $\BR$.
\end{prop}

The pairs in question arise because, loosely speaking, unlike
the $p$-adic case, where $-v_0=(p-1)v_{-1}+(p-1)v_{-2}+\dots$,
in the fibadic pattern we have two different representations:
$-u_0=u_{-1}+u_{-3}+u_{-5}+\dots=u_{1}+u_{-2}+u_{-4}+u_{-6}
+\dots$. The pairwise identification in question thus concerns
certain sequences that are finite to the left and cofinite to
the right. A formal way to establish this fact given in
\cite{Ver-ML} is as follows: while the standard representation
of the generators $u_n$ is $\pi(u_n)=\be^{-n}$, there is another
one, namely $\pi'(u_n)=(-\be)^n$. Then $\pi'(Q_\be)\subset\BR$,
and it suffices to show that every real number does have a
representation in base $((-\be)^{-n})_{n\in\BZ}$, and this
representation is unique everywhere except a certain countable
set. This claim follows from the results of Section~\ref{RE}
(see (\ref{N==})).

\begin{rmk} It is interesting to find out what will correspond
to different subsets of $\BR$ in $\wt X$. Since
$2=1+1=1+\be^{-1}+\be^{-2}=\be+\be^{-2}$, and similarly,
$3=\be^2+\be^{-2}$, etc., it is easy to see that $\BN\subset\wt X$
consists of finite sequences only.\footnote{``Finite" henceforward
will mean ``finite in both directions".} However, there are a lot
of finite sequences that do not yield a natural number, for
example, $1+\be^{-2}$. Moreover, it is shown in \cite{Ver-ML} that
if one takes the union of the finite sequences and the sequences
that are finite to the right and cofinite to the left, then after
the identification mentioned above, this set becomes naturally
isomorphic to the ring $\BZ[\be]\simeq \BZ+\BZ$ (and the finite
sequences are of course isomorphic to $\BZ[\be]\cap\BR_+$). This
is again the crucial difference with the $p$-adic case, where the
analogs are respectively $\BZ$ and $\BZ\cap\BR_+$.
\end{rmk}

\begin{rmk} More detailed results about the embedding of
different subsets of $\BR_+$ into $\wt X$ as well as about
relations with finite automata can be found in \cite{FrSa}. Note
also that by the theorem proven independently by A.~Bertrand
\cite{Ber77} and K.~Schmidt \cite{Sch80},
$\pi^{-1}(\BQ(\be)\cap\BR_+)$ is precisely the set of all
sequences finite to the left and eventually periodic to the
right (this is very similar to the $p$-adic case and is true for
all Pisot numbers).
\end{rmk}

The most important discovery made in \cite{Ver-ML} was the fact
that the fibadic analog of the solenoid $\mathcal S_p$ is
actually the 2-torus $\BT^2=\BR^2/\BZ^2$. Let us explain this in
detail as it appeared in subsequent works \cite{SV1, SV2}. Let
$\mathcal L_m$ stand for the Haar ($=$ Lebesgue) measure on
$\BT^m$, and $\Phi$ denote the {\em Fibonacci automorphism} of
$\BT^2$, namely, the algebraic automorphism given by the matrix
$$
M_\Phi=\begin{pmatrix}
  1 & 1 \\
  1 & 0
\end{pmatrix}.
$$
As is well known since the pioneering work by R.~L.~Adler and
B.~Weiss \cite{AdWe}, $\Phi$ is metrically isomorphic to the
two-sided $\be$-shift $\si_\be$. So, this is nothing new that
$\wt X$ as a set is essentially the torus; what {\bf is} new,
however, is that the natural arithmetic of $\wt X$ is the same
as the natural arithmetic of $\BT^2$. Our goal is thus dual: to
show that $\wt X$ is indeed {\em arithmetically} isomorphic to
the 2-torus (i.e., not only in the ergodic-theoretic sense but
in the arithmetic sense as well) and also to give a proof of the
Adler-Weiss Theorem cited above that reveals the arithmetic
structure of $\wt X_\be$. Both problems will be discussed
simultaneously.

We denote by $X_f$ the set of all sequences from $\wt X$ finite
to the left (recall that $X_f\simeq\BR_+$). Consider $x\ge0$ and
its greedy expansion $x=\sum_{k=-\infty}^\infty\e_k\be^{-k}$
given by (\ref{32-sided}) with $\e_k\equiv0$ for $k\le N(x)$.
Consider now the map $f_\be:X_f\to\BT^2$ acting by the formula
\[
f_\be(\e)=\{(\{x\},\{\be^{-1}x\})\mid x\ge0\},
\]
where $\{\cdot\}$ denotes the fractional part of a number. Let
$\mathcal R_\be\subset \BT^2$ denote the image of $X_f$ under
$f_\be$. Since $(1,\be^{-1})$ is an eigenvector of $M_\Phi$
corresponding to the eigenvalue $\be$, the set $W_\be$ is the
half-leaf of the unstable foliation for the Fibonacci
automorphism passing through $\0$. Hence
\begin{equation}
(f_\be\si_\be)(\e)=\Phi f_\be(\e)\label{dense}
\end{equation}
for any $\e=(\e_n)$ finite to the left.

Since the set $W_\be$ is dense in the 2-torus, as well as the
set of sequences finite to the right is dense in $\wt X$, we can
extend the relation~(\ref{dense}) to the whole compactum $\wt
X$, i.e. $(f_\be\si_\be)(\e)=\Phi f_\be(\e)$ everywhere on $\wt
X$. Besides, $f_\be$ is surjective and can be written in a very
``arithmetic" sort of way, namely
\begin{equation}
f_\be(\e)= \left(\sum_{k=-\infty}^{\infty}\e_k\be^{-k}\bmod 1,
\sum_{k=-\infty}^{\infty}\e_k\be^{-k-1}\bmod 1\right),\label{fbe}
\end{equation}
where the expression $\sum_{n=-\infty}^\infty x_n=x\bmod1$ means
that $\lim_{N}\|\sum_{n=-N}^Nx_n-x\|=0$. The number-theoretic
reason why these series do converge modulo 1 is that $\be$ is a
Pisot number, whence $\|\be^n\|\to0$ as an exponential rate.

\begin{lemma} \cite{SV1} The map $f_\be$ semiconjugates the
automorphisms $(\wt X, m, \si_\be)$ and $(\BT^2,\mathcal
L_2,\Phi)$, where $m$ denotes the (Markov) measure of maximal
entropy for $\si_\be$. Moreover, after an identification on $\wt
X$ that concerns a set of sequences of zero measure, $\wt X$
becomes an additive group $\wt X'$, and $f_\be$ becomes a group
homomorphism of $\wt X'$ and $\BT^2$.
\end{lemma}

Thus, we seemed to have succeeded in our attempt to insert the
arithmetic compactum $\wt X$ into $\BT^2$. However, this is not
that simple; the issue with $f_\be$ is that it is {\bf not}
bijective a.e. and thus cannot be regarded as an actual
isomorphism. In fact, in \cite{SV1} it was shown that it is
5-to-1 a.e.\footnote{This means that $\mathcal L_2$-a.e.
$x\in\BT^2$ has exactly 5 $f_\be$-preimages.} This is not a
coincidence -- in Section~\ref{PA} we will see that the
discriminant of an irrational in question plays an important
role in this theory (see Proposition~\ref{nnp}). The deep reason
why $f_\be$ has failed is because of the wrong choice of a
homoclinic point -- see below.

The way to construct an actual isomorphism is a slight
modification of $f_\be$. Namely, let $F_\be:\wt X\to\BT^2$ be
defined by the formula
\begin{equation}
F_\be(\e)=
\left(\sum_{k=-\infty}^{\infty}\e_k\frac{\be^{-k}}{\sqrt5}\bmod 1,
\sum_{k=-\infty}^{\infty}\e_k\frac{\be^{-k-1}}{\sqrt5}\bmod
1\right).\label{Fbe}
\end{equation}
Similarly to the above, the convergence of both series is a
consequence of the fact that $\|\be^n/\sqrt5\|=\be^{-n}/\sqrt5,\
n\ge0$.

\begin{thm}\cite{SV1} The map $F_\be$ is 1-to-1 a.e. It is both a
metric isomorphism of the automorphisms $(\wt X, m, \si_\be)$ and
$(\BT^2,\mathcal L_2,\Phi)$ and of the groups $\wt X'$ and
$\BT^2$.
\end{thm}

The question is, why $F_\be$ succeeded where $f_\be$ failed? The
reason becomes more transparent if we rewrite both maps. To do
so, we need to recall some basic notions and facts from
hyperbolic dynamics. Let $T$ be a hyperbolic automorphism of the
torus $\BT^m=\BR^m/\BZ^m$, $L_s$ and $L_u$ denote respectively
the leaves of the stable and unstable foliations passing through
$\0$. Recall that a point homoclinic to $\0$ (or simply a {\it
homoclinic point}) is a point which belongs to $L_s\cap L_u$. In
other words, $\t$ is homoclinic iff $T^n\t\to\0$ as
$n\to\pm\infty$. The homoclinic points are a group under
addition isomorphic to $\BZ^m$, and we will denote it by $H(T)$.
Each homoclinic point $\t$ can be obtained as follows: take some
$\mathbf n\in\BZ^m$ and project it onto $L_u$ along $L_s$ and
then onto $\BT$ by taking the fractional parts of all
coordinates of the vector (see \cite{Ver}).

We claim that both (\ref{fbe}) and (\ref{Fbe}) can be written in
the form
\begin{equation}
h_{\t}(\e)=\sum_{n\in\BZ}\e_nT^{-n}\t, \label{homo}
\end{equation}
where $T=\Phi$ and $\t=\t_1=(1,\be^{-1})$ in the case of $f_\be$
and $\t=\t_0=(1/\sqrt5,\be^{-1}/\sqrt5)$ in the case of $F_\be$.

The reason why $\t_0$ is ``better" than $\t_1$ is because it is
a {\em fundamental} homoclinic point, i.e., the one for which
the linear span of its orbit is the whole group $H(\Phi)$.

\begin{rmk}
Any fundamental homoclinic point for $\Phi$ is of the form
$\Phi^n\t_0$ for some $n\in\BZ$. In other words, $\t$ is
fundamental iff $\t=(\be^n/\sqrt5,\be^{n-1}/\sqrt5)\bmod\BZ^2$ for
some $n\in\BZ$. In the next subsection we will have a
generalization of this fact.
\end{rmk}

\subsection{Pisot automorphisms}\label{PA} The next step was
made by the author and A.~Vershik in \cite{SV2, VS2} -- it
concerned the general case of dimension~2. In this paper however
we will jump to the next stage, which will completely cover the
two-dimensional case, namely to the hyperbolic automorphisms of
the $m$-torus ($m\ge2$) whose stable (unstable) foliation is
one-dimensional.

Let $T$ be an algebraic automorphism of the torus $\BT^m$ given
by a matrix $M\in GL(m,\BZ)$ with the following property: the
characteristic polynomial for $M$ is irreducible over $\BQ$, and
a Pisot number $\be>1$ is one of its roots (we recall that an
algebraic integer is called \textit{a Pisot number}, if it is
greater than 1 and all its Galois conjugates are less than 1 in
modulus). Since $\det M=\pm1, \be$ is a \textit{unit}, i.e., an
invertible element of the ring $\BZ[\be]=\BZ[\be^{-1}]$. We will
call such an automorphism a \textit{Pisot automorphism}. Note
that since none of the eigenvalues of $M$ lies on the unit
circle, $T$ is hyperbolic. It is obvious that any hyperbolic
automorphism $T$ of $\BT^2$ or $\BT^3$ is either Pisot or one of
the automorphisms of the form $\pm T,\pm T^{-1}$ is such.

Our goal is, as above, to present a symbolic coding of $T$ which,
roughly speaking, reveals not just the structure of $T$ itself but
the natural arithmetic of the torus as well. Let us give a precise
definition.

\begin{Def} An \textit{arithmetic coding} $h$ of $T$ is a map from
$\wt X_\be$ onto $\BT^m$ that satisfies the following set of
properties:
\begin{enumerate}
\item $h$ is continuous and bounded-to-one;
\item $h\sigma_{\beta}=Th$;
\item $h(\e+\e')=h(\e)+h(\e')$ for any pair of sequences finite to
the left.
\end{enumerate}
\end{Def}

Thus, unlike the classical symbolic dynamics, where one has to
``encode" the action of $T$ itself, our goal is to give a {\bf
simultaneous} encoding of $T$ and the action of $\BT^m$ on
itself by addition. This makes the choice of $h$ much more
restricted; in fact, there are only a countable number of
arithmetic codings, as the following lemma shows:

\begin{lemma}\cite{SV2, S2} Any arithmetic coding of a Pisot
automorphism of $\BT^m$ is $h_\t$ given by (\ref{homo}), where
$\t\in H(T)$.
\end{lemma}

The issue is to find (if possible) an arithmetic coding of a
Pisot automorphism which is one-to-one a.e. We will call it a
{\em bijective arithmetic coding} or BAC. Before we formulate a
necessary and sufficient condition for $T$ to admit a BAC, we
need some auxiliary definitions. Let first the characteristic
equation for $\be$ be
\begin{equation*}
\beta^{m}=k_{1}\beta^{m-1}+k_{2}\beta^{m-2}+\cdots+k_{m},\quad
k_m=\pm1,
\end{equation*}
and $T_\be$ denote the toral automorphism given by the
\textit{companion matrix }$M_\be$ for $\be$, i.e.,
\[
M_\be=\left(
\begin{array}
[c]{ccccc}
k_1 & k_2 & \ldots &  k_{m-1} & k_m\\
1 & 0 & \ldots & 0 & 0\\
0 & 1 & \ldots & 0 & 0\\
\ldots & \ldots & \ldots & \ldots & \ldots\\
0 & 0 & \ldots & 1 & 0
\end{array}
\right).
\]
We need one more (arithmetic) condition on $\be$ to discuss. Let
$Fin(\be)$ denote the set of all $x\ge0$ having finite greedy
$\be$-expansion. It is obvious that
$Fin(\be)\subset\BZ[\be]_+=\BZ[\be]\cap\BR_+$. However, the
inverse inclusion does not holds for some Pisot units; those for
which it does hold, are called {\em finitary}. For examples see,
{\em e.g.}, \cite[\S2]{S1}. The property of $\be$ to be finitary
helps in many Pisot-related issues, but our goal here is to
present a more general result, which is based on a more general
property.
\begin{Def}
A Pisot unit $\be$ is called \textit{weakly finitary} if for any
$\de>0$ and any $x\in\mathbb{Z}[\be]_+$ there exists $f\in
Fin(\be)\cap(0,\de)$ such that $x+f\in Fin(\beta)$ as well.
\label{wf}
\end{Def}

This notion has appeared in different contexts and is related to
different problems -- see \cite{Ak,Hol,S2}. The following
conjecture (apparently, very difficult to prove) is shared by
most experts.

\begin{conjecture}
\label{conject}Any Pisot unit is weakly finitary.
\end{conjecture}

To find out more about this property and about the algorithm how
to verify that a {\bf given} Pisot unit is weakly finitary, see
\cite{Ak}.

Return to our setting. We assume the following conditions to be
satisfied:
\begin{enumerate}
\item $T$ is algebraically conjugate to $T_\be$, i.e.,
there exists a matrix $C\in GL(m,\mathbb{Z})$ such that
$CM=M_\be C$ (notation: $T\sim T_\be$).
\item A homoclinic point $\t$ is fundamental.
\item $\be$ is weakly finitary.
\end{enumerate}

\begin{thm}\label{Pisenc} \textrm{(1)} If a Pisot automorphism
$T$ admits a BAC, then $T$ is algebraically conjugate to
$T_\be$.\newline \textrm{(2)} Assume that the three conditions
above are satisfied. Then $T$ admits an arithmetic coding
bijective a.e.
\end{thm}

\begin{rmk} Theorem~\ref{Pisenc}~(2) for the case of
finitary Pisot eigenvalue has been proven by Le~Borgne in his
Ph.~D. Thesis \cite{Leb} (see also \cite{Sch} for some cases).
\end{rmk}

\begin{rmk} If $T\sim T_\be$, then a fundamental homoclinic point
always exists. Thus, modulo Conjecture~\ref{conject}, the
algebraic conjugacy to the companion matrix is the necessary and
sufficient condition for a Pisot automorphism to admit a BAC.
\end{rmk}

For the rest of the subsection we assume $\be$ to be weakly
finitary. Similarly to the Fibonacci case, the set $\wt X_\be$
is an {\em almost group} in the following sense.

\begin{prop}\cite{S2}
\label{argroup} Let $\sim$ denotes the identification on $\wt
X_\be$ defined as follows: $\e\sim\e'$ iff $h_\t(\e)=h_\t(\e')$,
where $\t$ is fundamental. Then it touches only a set of measure
zero, and $\wt X_\be'=\wt X_\be/\sim$ is a group isomorphic to
$\BT^m$.
\end{prop}

A natural question to ask is as follows: what is the number of
preimages of a generic point if $\t$ is not fundamental? (for
instance, if $T\not\sim T_\be$) In \cite{S2} this question is
answered completely.

We start with the case $T=T_\be$ and show how this problem is
related to Algebraic Number Theory. Let
$$
\P_\be=\{\xi\in\BR: \|\xi\be^n\|\to0,\,\, n\to+\infty\}.
$$
It is well-known that $\P_\be\subset\BQ(\be)$ (see, {\em e.g.},
\cite{Cas}). Let $\Tr(\xi)$ denote the trace of $\xi$, i.e. the
sum of $\xi$ and all its conjugates. It is shown in \cite{S1}
that the set $\P_\be$ is a commutative group under addition
containing $\BZ[\be]$ and also that it can be characterized as
follows:
$$
\P_\be=\{\xi\in\BQ(\be): \Tr(a\xi)\in\BZ\ \text{\rm{for any}}\
a\in\BZ[\be]\}.
$$

\begin{lemma}\cite{S1} There exists a one-to-one correspondence
between the homoclinic points and the elements of $\P_\be$.
Namely, $\t\in H(T)$ if and only if
$$
\t=(\xi,\xi\be^{-1},\dots, \xi\be^{-m+1})\bmod\BZ^m
$$
for some $\xi\in\P_\be$.
\end{lemma}

Thus, any arithmetic coding of $T_\be$ is of the form
\begin{equation}
h_\xi(\e)=\sum_{k\in\BZ}\e_kT^{-k}\t=
\lim_{N\to+\infty}\left(\sum_{k=-N}^\infty
\e_k\be^{-k}\right)\begin{pmatrix} \xi\\ \xi\be^{-1}\\ \vdots\\
\xi\be^{-m+1}\end{pmatrix},
\end{equation}
where $\xi=\xi(\t)\in\P_\be$. Let $N(\cdot)$ denote the norm in
$\BQ(\be)$ and $D=D(\be)$ stand for the discriminant of $\be$.

\begin{prop}\label{nnp}\cite{S2} The map $h_\xi$ is $K$-to-1
a.e., where $K=|DN(\xi)|$.
\end{prop}

\begin{rmk} Thus, $h_\xi$ is a BAC if and only if
$N(\xi)=\pm 1/D$, which is equivalent to the fact that
$\xi/\xi_0$ is a unit in $\BQ(\be)$. If $\xi=1$, then we come to
the historically the first attempt to encode a Pisot
automorphism undertaken by A.~Bertrand-Mathis in \cite{Ber1}.
Now we see that $h_1$ is in fact $|D|$-to-1 (provided $\be$ is
weakly finitary).
\end{rmk}

Consider now the general case. We will be interested in the {\bf
minimal} number of preimages of $h_\t$ that one can attain for a
given $T$. Let $M\in GL(m,\BZ)$ denote the matrix which
determines $T$. To answer the above question, we are going to
describe {\bf all} integral square matrices that semiconjugate
$M$ and $M_\be$. Let for $\mathbf{n}\in\BZ^m$ the matrix
$B_M(\mathbf{n})$ be defined as follows (we write it
column-wise):
\[
\begin{aligned}
B_M(\mathbf{n})=(&M\mathbf{n},(M^2-k_1 M)\mathbf{n}, (M^3-k_1
M^2-k_2 M)\mathbf{n},\ldots,\\&M^{m-1}-k_1
M^{m-2}-\cdots-k_{m-2}M)\mathbf{n}, k_m\mathbf{n}).
\end{aligned}
\]

\begin{lemma}\cite{S2}
\label{Semic}Any integral square matrix satisfying the relation
\begin{equation*}
BM_{\beta}=MB
\end{equation*}
is $B=B_M(\mathbf{n)}$ for some $\mathbf{n}\in\BZ^m$.
\end{lemma}

Let
\[
f_{M}(\mathbf{n}):=\det B_{M}(\mathbf{n})
\]
(an $m$-form of $m$ variables).

\begin{prop}
\label{numpreimage} Let $\t\in H(T)$. Then there exists
$\mathbf{n}\in\BZ^m$ such that
\[
\#\varphi_{\mathbf{t}}^{-1}(x)\equiv|f_{M}(\mathbf{n})|
\]
for $\mathcal{L}_m$-a.e. point $x\in\BT^m$.
\end{prop}

\begin{cor} The minimal number of preimages for an arithmetic
coding of $T$ equals the arithmetic minimum of the form $f_M$.
\end{cor}

Thus, $T$ admits a BAC iff the Diophantine equation
\begin{equation}\label{pmm1}
f_M(\mathbf n)=\pm1
\end{equation}
is solvable. In the case $m=2$, (\ref{pmm1}) is especially
natural: if $M=\begin{pmatrix} a & b \\ c & d\end{pmatrix}$,
then it is
\[
cx^2-(a-d)xy-by^2=\pm1
\]
and therefore, belongs to the class of well-known quadratic
Diophantine equations. For more details about the
two-dimensional case see \cite{SV2, VS2}.

\subsection{General case} The previous subsection has covered
the case when one of the eigenvalues of the matrix of an
automorphism is outside (inside) the unit disc and all the
others are inside (resp. outside). The model explained above
looks rather natural, explicit and canonical. What can be done
in case when at least two eigenvalues are outside the unit disc
and at least two -- inside it? The main difficulty here lies in
the fact that unlike the Pisot case, where the entropy is
$\log\be$ and the $\be$-compactum is the obvious candidate for a
coding space, in the general case this choice is not at all
obvious.

There are several constructions that cover the general
hyperbolic (or even ergodic) case, and each of them has its own
advantages and disadvantages. Before we describe all of them in
detail, let us try to understand what is that we actually want
from an arithmetic coding. Obviously, there are no new
properties of algebraic toral automorphisms that can be revealed
this way -- simply because they all are so well
known.\footnote{For instance, the construction of Markov
partitions for the hyperbolic automorphisms of a torus (even for
more general Axiom A diffeomorphisms \cite{Sinai, Bowen}) was
revolutionary in the sense that although it was practically
implicit, it nonetheless allowed to show ``for free" (with the
help of the famous Ornstein Theorem, of course) that they are
all Bernoulli, which completely justified all the hard efforts
and technicalities.} What then? The unclear situation with this
has, in my opinion, led to a certain impasse in this theory. No
model seems to be canonical, and until we find an appropriate
application, any theory will be a
$\,\mathfrak{Ding}$-$\mathfrak{an}$-$\mathfrak{sich}$.

Let us also note that there are two main challenges any general
arithmetic encoding has to meet:
\begin{enumerate}
\item it has to be bounded-to-one and, if possible, one-to-one
a.e.;
\item the alphabet - it should be as simple as possible
(preferably integers).
\end{enumerate}

Which one is more important (if one cannot achieve both aims)?
Here is one possible application that might measure the value of
different constructions.

We have already mentioned the theorem on maps with holes proven
by S.~Bundfuss, T.~Krueger and S.~Troubetz\-koy in \cite{BKT}
(see Section~\ref{un}). Recall that this theorem claims that a
if one cuts out a ``typical" parallelepiped from $\BT^m$ along
the directions of the stable and unstable foliations with a
vertex at $\0$ and the sides of length $a_1,\dots,a_m$, then the
corresponding exclusion map will be a subshift of finite type.
This nice result however does not give any conditions on $a_i$
for this subshift to be {\bf nondegenerate}. At the same time,
it is possible to show that if $a_i$ are very small, then its
entropy will be positive, and it is obvious that for ``large"
$a_i$ the images of the hole will cover the whole torus, so it
will be degenerate. Thus, if we make a natural assumption that
similarly to the one-dimensional case, the entropy of the
exclusion subshift is a continuous function of
$(a_1,\dots,a_m)$, then there exists a threshold similar to the
Komornik-Loreti constant for the map $T_\be$ (see
Section~\ref{un}). In other words, we will have the surface
$\Pi$ in the space $(a_1,\dots,a_m)$, underneath which the
entropy of the exclusion subshift parameterized by $(a_i)$ is
positive, and it is zero above $\Pi$.

We do not know how the surface $\Pi=\Pi(T)$ looks like even in
the case of the Fibonacci automorphism (where it in fact must be
a curve). Nonetheless, we believe the exact simple formula for
the symbolic encoding like (\ref{homo}) with an explicitly
described symbolic compactum will probably help to reformulate
the problem in terms symbolic sequences and to treat it in a way
similar to the one described in \cite{GS1}. In particular, let
us ask the following question: is there any multidimensional
analog of the Thue-Morse sequence (cf. Section~\ref{BE})?

For this problem it is obvious that a bounded-to-one encoding
map will be sufficient, as long as the set of digits and the map
itself are explicit (because the entropy is preserved). We plan
to return to this problem in our subsequent papers. Now it is
time to present all the models known to date and to compare
them.

\subsubsection{The construction of Kenyon and Vershik}
Historically the first general arith\-metic symbolic model for
the hyperbolic automorphisms was suggested by R.~Kenyon and
A.~Vershik \cite{KenVer} (published in 1998 but written in
1995). This model is based on certain constructions that
intensively use Algebraic Number Theory. We refer the reader to
the textbooks, {\em e.g.}, \cite{BorShaf, FrolTa} for the
relevant notions and results. We will keep the original notation
of \cite{KenVer} and hope this will not make any confusion with
the notation of the rest of the present paper.

\smallskip\noindent\texttt{Alphabet.} Let $\la_1,\dots,\la_m$
be the eigenvalues of $M$, where $|\la_i|>1$ if and only if
$i=1,\dots,k$. Let $K=\BQ[x]/p(x)$, where $p(x)$ is the
characteristic polynomial for $M$, and $\mathcal O$ denote the
ring of integers in $K$. The ring $K$ (and therefore, $\mathcal
O$ as well) is naturally embedded into $\BR^m$ via the standard
coordinate-wise embeddings. The set $\mathcal O$ becomes a
full-rank lattice in $\BR^m$.

We denote by $B$ the closed ball centered at $\0$ with the
radius $r$ defined as the smallest $t$ such that its any
translation has a nonempty intersection with $\mathcal O$.
Finally, $D:=\mathcal O\cap(B+xB)$, where multiplication by $x$
symbolizes the multiplication by the companion matrix for $M$.
The set $D$ is shown to be finite, and this is precisely the set
of digits for the model of \cite{KenVer}.

\smallskip\noindent\texttt{Coding.} Let $\si$ denote the shift
on $D^{\BN}$ and $(\Si_u,\si)$ denote the subshift defined as
follows: assume $D$ is endowed with some full order $\prec$;
this creates the lexicographic ordering on $D^{\BN}$. If
$(\e_1,\dots,\e_j)$ is a finite sequence, we say it is
non-minimal if there exists a word
$(\e'_1,\dots,\e'_j)\prec(\e_1,\dots,\e_j)$ such that
$\sum_{i=1}^j\e_ix^{j-i}=\sum_{i=1}^j\e'_ix^{j-i}$. If a
sequence is not non-minimal, we call it {\em minimal}.

The space $\Si_u$ is thus the closed shift-invariant subset of
$D^{\BN}$ consisting of those sequences whose finite
subsequences are all minimal. The coding space will be
$(\Si,\si)$, where $\Si$ is the natural extension of $\Si_u$.

\begin{prop}\cite{KenVer} The subshift $(\Si,\si)$ is sofic.
\end{prop}

Now let us follow the authors of \cite{KenVer} in their
construction of the encoding map. Define for
$d=(d_0,d_1,\dots)\in D^{\BN}$,
\[
S_i(d)=\sum_{j=0}^\infty\rho_i(d_j)\la_i^{-j}
\]
if $i=1,\dots,k$ and
\[
S_i(d)=\sum_{j=0}^\infty\rho_i(d_j)\la_i^{j+1}
\]
otherwise. Furthermore, let $R_u:D^{\BN}\to W_u$ (the unstable
eigenspace of the companion matrix) act as follows:
$R_u(d)=(S_1(d),\dots,S_k(d))$ and similarly
$R_s(d)=(S_{k+1}(d),\dots,S_m(d))$. Finally, let
\[
R(\dots,d_{-1},d_0,d_1,\dots):=
R_u(d_0,d_1,\dots)-R_s(d_{-1},d_{-2},\dots)
\]
be the map from $\Si$ to $\BR^m$, and $\pi$ denote the natural
projection from $\BR^m$ to $\BT^m$.

\begin{thm}\cite{KenVer} The map $\pi R$ is a factor map from
$(\Si,\si)$ to $(\BT^m,T)$. It is bounded-to-one everywhere and
constant-to-one a.e.
\end{thm}

\noindent\texttt{Examples.} The authors consider in detail the
Fibonacci and similar quadratic cases as well as some cubic
cases. Unfortunately, none of them uses the original set of
digits $D$ described above (in the Fibonacci case, for example,
they take the conventional $D=\{0,1\}$). Thus, it is difficult
to assess the effectiveness of this model; nonetheless, the
authors show how to deal with the ``reasonable" choice of digits
in specific cases. Note also that E.~Hirsch proved in
\cite{Hirsch} that it is impossible for a general case to use
this model with $D$ containing just nonnegative integers.

\subsubsection{The construction of Le Borgne} The model
suggested by S.~Le~Borgne in his Ph.~D.~Thesis \cite{Leb} (see
also \cite{Leb1, Leb2}) is in fact a generalization (map-wise)
of the Pisot model described above. As usual, we preserve the
author's notation.

\smallskip\noindent\texttt{Alphabet.} Let $F_u, F_s$ denote the
unstable and stable foliations for $T$ and $\pi_u$ stand for the
projection from $\BR^m$ onto $F_u$ along $F_s$, and we define
$\pi_s$ in a similar way. Let $M_u$ denote the restriction of
$M$ to $F_u$.

Assume $E\subset\pi_u(\BZ^m)$ to be a finite set, and
\begin{equation}
W_E=\left\{\sum_{j=0}^\infty M_u^{-j}e_j\mid e_j\in
E\right\}.\label{Mu}
\end{equation}

\begin{lemma}\label{cen}\cite{Leb, Leb1} It is always possible
to choose $E$ in such a way that the interior of $W_E$ is
nonempty.
\end{lemma}

Henceforward we assume $E$ to be such, and $W=W_E$. Let now $Y$
denote the set of all sequences that appear in the
expansion~(\ref{Mu}) and let $Z$ be its natural extension.
Finally, denote by $X$ the maximal transitive subshift of $Z$.

\begin{lemma}\cite{Leb, Leb1} The shift $(X,\si)$ is sofic and
has a unique measure of maximal entropy ($\nu$, say).
\end{lemma}

The set $X$ is the sought symbolic compactum. The ``digits" thus
are in fact vectors, and the actual choice is hidden in
Lemma~\ref{cen}; see below how to convert vectors into (more
conventional) integers in the case of $M$ algebraically
conjugate to its companion matrix.

\smallskip\noindent\texttt{Coding.} Let $X$ be as above, and
$\phi:X\to\BT^m$ be defined by the formula
\begin{equation}
\phi(\e)=\sum_{j=1}^\infty\pi_u(M^{-j}\e_j)\bmod\BZ^m-
\sum_{j=-\infty}^0\pi_s(M^{-j}\e_j)\bmod\BZ^m.\label{pphi}
\end{equation}

\begin{thm}\cite{Leb, Leb1} The map $\phi$ given by (\ref{pphi})
is surjective, H\"older continuous and $p$-to-one a.e for a
certain $p\in\BN$. It semiconjugates the transitive sofic shift
$(X,\nu,\si)$ and $(\BT^m,\mathcal L_m,T)$.
\end{thm}

The main issue is to make it one-to-one a.e. (by an appropriate
choice of $E$) as well as to make the alphabet more canonical.
In the case when $M$ is algebraically conjugate to its companion
matrix, this has been partially done in the thesis \cite{Leb}.
Let $\Xi=\pi_u^{-1}(E)\subset \BZ^m$.

\begin{prop} Let $u_0\in\BZ^m$ be such that $\langle M^ju_0\mid
j=0,1\dots,m-1\rangle=\BZ^m$. There exists $N\ge1$ such that
$\Xi$ may be chosen in the form $\{-Nu_0,\dots,Nu_0\}$.
\end{prop}

Thus, in a way, one might say that the digits are integers. The
author also shows how (theoretically) the alphabet can be
constructed but gives no non-Pisot examples.

\subsubsection{The construction of Schmidt} The paper \cite{Sch}
by K.~Schmidt appeared right after \cite{SV2} and used the map
defined by (\ref{homo}). More precisely, the case considered in
\cite{Sch} was more general than the hyperbolic toral
automorphisms: the author deals with expansive group
automorphisms of compact abelian groups. We will not be
concerned with the general case though and will confine
ourselves to the setting in question.

\begin{thm}\cite{Sch} For a given hyperbolic automorphism
$T$ of $\BT^m$ whose matrix is algebraically conjugate to its
companion matrix there exists a topologically mixing sofic
subshift $V$ of $l^\infty(\BZ,\BZ)$ such that
\begin{enumerate}
\item $h_{\t}(V)=\BT^m$, where $h_\t$ is given by (\ref{homo});
\item The restriction of $h_\t$ to $V$ is one-to-one everywhere
except the set of doubly transitive points of $T$.
\end{enumerate}
\end{thm}

\begin{rmk} The proof given in \cite{Sch} is non-constructive.
As the author himself states, the above theorem only asserts the
{\em existence} of a sofic shift $V$ with the properties
described above.
\end{rmk}

\subsubsection{Conclusions} Let us compare all models by
gathering all we know about them in the following table:

\bigskip
\noindent
\begin{tabular}{|p{2.9truecm}|p{2.9truecm}|p{2truecm}|
p{2truecm} |p{2truecm}|} \hline \hphantom{XXXX} &
\textbf{Sidorov-Vershik} & \textbf{Kenyon-Vershik} & \textbf{Le
Borgne} & \textbf{Schmidt} \\
\hline\hline \textit{Auto\-mor\-phisms covered}& 2D and
generalized Pisot (modulo arith\-me\-tic con\-jecture) &
Hyperbo\-lic & Hyperbo\-lic & Hy\-per\-bo\-lic cyclic\\
\hline \textit{Is the sub\-shift explicit?}
& Yes & No & No & No \\
\hline \textit{Is the coding canonical?} & Yes & Yes & No & No
\\
\hline \textit{``Digits"} & Nonnega\-tive integers & Algebra\-ic
numbers & Vectors & Integers \\
\hline \textit{The encoding map is} & $K$-to-one a.e. and
one-to-one a.e. for the cyclic & $K$-to-one a.e. & $K$-to-one
a.e. & Boun\-ded-to-one and one-to-one a.e.\\
\hline
\end{tabular}

\begin{rmk} Here a {\em generalized Pisot automorphism} means
that its stable (unstable) foliation is one-dimensional. They
all can be arithmetically encoded using the construction for the
Pisot automorphisms -- see \cite{S2} for details. The expression
``$K$-to-one a.e." implies that there exists $K\in\BN$ such that
almost every point of the torus has $K$ preimages and ``cyclic"
means ``the matrix is algebraically conjugate to its companion
matrix".
\end{rmk}

It is also worth noting that an attempt to deal with the general
case has been undertaken by the author in \cite{S2}. The idea is
as follows: assume $S$ is a hyperbolic automorphism of $\BT^m$
and $T$ is a generalized Pisot automorphism of $\BT^m$ that
commutes with $S$. Then $S=\sum_{j=0}^{m-1}c_jT^j$ with
$c_j\in\BQ$. Actually, the denominators of $c_j$ are known to be
bounded, and we assume that $c_j\in\BZ$ for all $j$. Recall that
the map $h_\t$ given by (\ref{homo}) semiconjugates (or
conjugates if $T$ is cyclic) the shift $(\wt X_\be,\si_\be)$ and
$(\BT^m,T)$. Hence the same map semiconjugates the linear
combination of the powers of $\si_\be$, namely,
$\sum_{j=0}^{m-1}c_j\si_\be^j$, and $S$ (recall that by
Proposition~\ref{argroup} the set $\wt X_\be$ is an ``almost
group", whence any fixed finite integral combination of the
powers of the shift is well defined a.e.). Thus, if we do not
require that it must be necessarily a shift that encodes $S$, we
are practically done.

The main issue is number-theoretic: the question is whether in a
given algebraic field $K\simeq\BQ(\la)$ there exists a Pisot
unit $\be$, and if it exists, whether it can be found in such a
way that $\la$ is an {\bf integral} linear combination of powers
of $\be$. Of course, if, for instance, $K$ is totally real
(which leads to the {\em Cartan action}, i.e., the
$\BZ^{m-1}$-action by algebraic automorphisms), then it always
contains a Pisot unit but the second property seems to be more
difficult to prove -- it requires some knowledge about the
structure of the Pisot units in an algebraic field, which is
apparently missing in the classical Algebraic Number Theory.

\begin{exam}\cite{S2} Let $M=\left(
\begin{array}
[c]{cccc}
4 & 0 & -3 & 1\\
1 & 0 & 0 & 0\\
0 & 1 & 0 & 0\\
0 & 0 & 1 & 0
\end{array}
\right)$. Note $M$ is a companion matrix, and its spectrum is
purely real. Now take the action generated by $M_1=M, M_2=M+E$
and $M_3=M-E$. It is easy to check that they all belong to
$GL(4,\mathbb{Z})$ and that this will yield a Cartan action on
$\mathbb{T}^4$ as well as the fact that the dominant eigenvalue
$\beta$ of $M$ is indeed weakly finitary. We leave the details
to the reader. Therefore, the usual mapping $h_\t$ conjugates
the action generated by $(\sigma_{\beta},\sigma_{\beta}+id,
\sigma_{\beta}-id)$ on the compactum $X_{\beta}$ and the Cartan
action generated by $(T_1,T_2,T_3)$. Furthermore, $T_3$ has two
eigenvalues strictly inside the unit disc and two strictly
outside it. Perhaps, this is the first ever explicit bijective
a.e. encoding of a non-generalized Pisot automorphism (though
not by means of a shift).
\end{exam}

Is this model any good application-wise? I am not sure; in
particular, for the maps with holes the fact that instead of a
shift we have this modified map, does not help a lot. However,
it might be worth trying to apply it, when the Pisot case
becomes clear. The author is grateful to A.~Manning,
M.~Einsiedler and K.~Schmidt for helpful discussions and
number-theoretic insights regarding this question.

\section*{Appendix: adic transformations}

In this appendix we are going to describe the class of maps on
symbolic spaces which is in a way transversal to the
shifts.\footnote{Actually, this statement can be made precise
whenever the symbolic space is stationary ($=$ shift-invariant)
-- see, {\em e.g.}, \cite{Ver93} for some cases. As we will see,
the adic transformations cover a much wider class of spaces.}
Let us give the precise definition.

Let $(D_k)_{k=1}^\infty$ be a sequence of finite sets, $r_k=\#
D_k$, and let $\X':=\prod_1^\infty D_k$ endowed with the weak
topology. A closed subset $\X$ of $\X'$ is called a {\em Markov
compactum} if there exists a sequence of 0-1 matrices
$(M^{(k)})_{k=1}^\infty$, where $M^{(k)}$ is an $r_k\times
r_{k+1}$ matrix, such that
\[
\X=\X(\{M^{(k)}\})=\{(x_1,x_2,\dots)\in\X':M^{(k)}_{x_kx_{k+1}}
=1\}.
\]
In other words, $\X$ is a (generally speaking, non-stationary)
analog of topological Markov chain, and the $M^{(k)}$ are its
incidence matrices. Assume that there is a full ordering $\prec_k$
on each set $D_k$. Then this sequence of orderings induces the
partial lexicographic order on $\X$ in a standard way: two
distinct sequences $x$ and $x'$ are comparable iff there exists
$n\ge1$ such that $x_n\neq x_n'$ and $x_k=x'_k$ for all $k\ge
n+1$. Then $x\prec x'$ iff $x_n\prec_n x'_n$.

\begin{Def} The {\em adic transformation} $S$ on $\X$ is defined
as a map that assigns to a sequence $x$ its immediate successor in
the sense of the lexicographic ordering defined above (if exists).
\end{Def}

\begin{rmk} If $M^{(k)}_{x_kx_{k+1}}=1$ for {\bf all} pairs
$(x_k,x_{k+1})$, then we have the {\em full odometer} or the {\em
$\mathbf{r}$-adic transformation} $S_{\mathbf r}$. If we identify
$D_k$ with $\{0,1,\dots,r_k-1\}$, then $S_{\mathbf r}$ acts on the
the set of $\mathbf r$-adic integers in the following way: every
finite sequence $x$ can be associated with a nonnegative integer
as usual, i.e., $N=\sum_k x_k r_k$; then $S_{\mathbf r}(N)=N+1$.
The map $S_{\mathbf r}$ on the whole space is thus the profinite
completion of the operation $N\mapsto N+1$. It is well defined
everywhere except the sequence $(r_1-1,r_2-1,\dots)$.

As is well known, this map has purely discrete spectrum for any
$\mathbf{r}=(r_k)_1^\infty$. Thus, the adic transformation on an
arbitrary Markov compactum may be regarded as a Poincar\'e map for
the full odometer and some Markov subcompactum.
\end{rmk}

The adic transformation is known to be well defined a.e. for a
large class of systems (see, {\em e.g.}, \cite{LivVer}). The
importance of this model is confirmed by the following theorem
proved by A.~Vershik in the seminal paper \cite{Ver82}, where
the notion in question was first introduced (see also
\cite{Ver81}).

\begin{thm} \cite{Ver82} Each ergodic automorphism of the Lebesgue
space is metrically isomorphic to some adic transformation.
\label{metric_iso}
\end{thm}

\begin{rmk}
A ``topological" version of this theorem have been obtained by
M.~Herman, I.~Putnam and C.~Skau \cite{HPS}. Recently A.~Dooley
and T.~Hamachi obtained a version of this theorem for the
quasi-invariant measures of type~III \cite{DoHam}. Note also
that the adic realization in a special case were earlier
considered by M.~Pimsner and D.~Voiculescu \cite{PimVo} in
connection with approximations of certain operator algebras.
\end{rmk}

Note that although the proof of Theorem~\ref{metric_iso} is
based on Rokhlin's Lemma and is thus to some extent
constructive, there are very few explicit examples of ``adic
realization". The irrational rotations of the circle are among
those rare exceptions (see Section~\ref{RE}); unfortunately,
even for a general ergodic shift on the 2-torus the model seems
to be hardly constructible.

One more fact worth noting is that if a Markov compactum is
stationary (i.e., if $M^{(k)}\equiv M$ for any $k\ge1$), then,
as was shown by Livshits \cite{Liv}, the adic transformation on
it is isomorphic to a {\em substitution} or, as it is more
appropriate to call it, a {\em substitutional dynamical system}.
The converse is also true, i.e., any primitive substitution has
a stationary adic realization. For instance, the Fibonacci
substitution $0\to01,\ 1\to0$ is isomorphic to the adic
transformation on $X_G$, while the Morse substitution $0\to01,\
1\to10$ leads to the the adic transformation on
$\prod_1^\infty\{0,1\}$ with the alternating ordering similar to
the one described in Section~\ref{RE-1} (see the second model).
A good exposition of this theory can be found in \cite{Ver93}.

\end{document}